\documentclass[review]{elsarticle}

\usepackage{lineno,hyperref}
\usepackage{amsmath}
\usepackage{graphicx}
\usepackage{subfigure}
\usepackage{caption}
\usepackage{float}
\usepackage{amsfonts}
\usepackage{amssymb}
\usepackage{indentfirst}
\usepackage{booktabs}
\usepackage{multirow}
\usepackage{color}
\usepackage{epstopdf}
\modulolinenumbers[5]

\usepackage[titletoc,title]{appendix}
\newtheorem{definition}{Definition}[section]

\journal{Journal of \LaTeX\ Templates}

\biboptions{authoryear}

\begin{document}
%
%
\def\ear{\textit{E\/}}       
\def\lat{\textit{L\/}}       
\def\dur{\textit{D\/}}       
\def\wght{\textit{w}}        

\def\usg{\textit{A\/}}       
\def\delt{\Delta\/}          

\def\TW{\textit{TW}}         
\def\tw{\textit{tw}}         
\def\RTW{\textit{RTW}}       
\def\rn{\textit{rn}}         
\def\rtw{\textit{rtw}}       
\def\srn{\textit{srn}}       
\def\srtw{\textit{srtw}}     
\def\WBeg{\textit{Beg\/}}    
\def\WEnd{\textit{End\/}}    

\def\SBeg{S_{\textit{Beg}}}  
\def\SEnd{S_{\textit{End}}}  

\def\conf{\textit{conf}}     
\def\paon{\textit{paon}}     
\def\paot{\textit{paot}}     

\def\mVC{\textit{mVC}}       
\def\mVB{\textit{mVB}}       
\def\mC{\textit{mC}}         

%
%
\def\OBeg{\textit{t}}        

\begin{frontmatter}

\title{A Mixed Integer Linear Programming Model\\
       for Multi-Satellite Scheduling}

\author[mymainaddress,mysecondaryaddress]{Xiaoyu Chen}
\ead{xiaoyu.chen@cug.edu.cn}

\author[mysecondaryaddress]{Gerhard Reinelt}
\ead{ip121@uni-heidelberg.de}

\author[mymainaddress]{Guangming Dai\corref{mycorrespondingauthor}}
\cortext[mycorrespondingauthor]{Corresponding author}
\ead{cugdgm@126.com}

\author[mysecondaryaddress]{Andreas Spitz}
\ead{spitz@informatik.uni-heidelberg.de}

\address[mymainaddress]{School of Computer Science, China University of Geosciences, Wuhan 430074, China}
\address[mysecondaryaddress]{Institute of Computer Science, Heidelberg University, Heidelberg 69120, Germany}

\begin{abstract} 
We address the \emph{multi-satellite scheduling problem with limited observation capacities} that arises from the need to observe a set of targets on the Earth's surface using imaging resources installed on a set of satellites. We define and analyze the conflict indicators of all available visible time windows of missions, as well as the feasible time intervals of resources. The problem is then formulated as a mixed integer linear programming model, in which constraints are derived from a careful analysis of the interdependency between feasible time intervals that are eligible for observations. We apply the proposed model to several different problem instances that reflect real-world situations. The computational results verify that our approach is effective for obtaining optimum solutions or solutions with a very good quality.
\end{abstract}

\begin{keyword}
Scheduling\sep Earth Observing Satellites\sep Integer Programming\sep Mathematical Programming
\end{keyword}

\end{frontmatter}


\section{Introduction}\label{Sec:Introduction}
Earth-observing satellites~(EOS) are specially designed for the observation of activities or areas on the Earth's surface, and play an increasingly important role in resource explorations, disaster alerts, environmental damage analysis, and many other imaging demands~\citep{LiuX2017}. An EOS can photograph the target with a variety of equipped resources, such as sensors or cameras. Each resource has a limited observation region on the Earth's surface that is formed by the subpoint of the satellite's resource and the field of view. The observation activity can be controlled by the swing angle and the rotation angle of the resource~\citep{HabetD2010,LiuX2014}. Clearly, it is only possible for the resource to accomplish the observation if the target is visible to it~\citep{MaoT2012}. 

During the observation process, every target has to be observed for a specified duration that depends on the resource, which can be calculated from the orbiting speed of the satellite and the scanning speed of the resource~\citep{NiuX2015}. The observation operation must be continuously and completely executed within a time window during which the target is visible to the satellite~\citep{YaoF2010}. For each mission, there may exist multiple feasible observation windows per resource. Furthermore, some additional constraints may need to be taken into account, such as operational constraints of satellites, energy capacity restrictions, resource availability, or requirements of special resource types. Additionally, the swing angle and rotation angle of the resource must be set to point at the target. Thus, a setup time between two consecutive successful observations has to be considered to adjust the orientation of the resource~\citep{MaoT2012}. While the number of EOS is continuously increasing, so is the number of observation requests. Therefore, given the cost of operating satellites, it is reasonable to assume that the capacity of satellites to satisfy customer demands for observation missions is a scarce resource~\citep{WuG2012}, and that it may not be possible to satisfy all mission demands during a given observation period. Thus, the development of effective scheduling approaches is pertinent, such as the approach we discuss in the following.

In the following, we address the multi-satellite scheduling problem with limited observation capacities. In comparison to related work, we make four primary contributions. (1) We analyze the capacity of the resources and the distribution of visible time windows of missions. (2) We introduce and define conflict indicators of available visible time windows of missions. (3) We derive further constraints from a careful analysis of the interdependency of time intervals that are eligible for observations. (4) Finally, we formulate the problem as a mixed integer linear programming model. Our computational results indicate the effectiveness and efficiency of the proposed method.

The remainder of this paper is organized as follows. We first review the related work and the state-of-the-art in multi-satellite scheduling in Section~\ref{Sec:RelatedWork}. By defining the conflict indicator of all available visible time windows of missions and by analyzing the capability of all feasible time intervals of resources, we formulate an exact mixed-integer linear program~(MILP) in Section~\ref{Sec:Model}. Simulation results and a performance analysis on a series of benchmark problem instances are given in Section~\ref{Sec:Simulation}. Finally, a summary and our conclusions are provided in Section~\ref{Sec:Conclusion}. 

\section{Related Work}\label{Sec:RelatedWork}
Given the complexity of the issue, a large portion of previous works is concerned with single satellite scheduling and address the efficient performance by providing an optimal solution and an upper bound. A common set of benchmark instances~(S5-DPSP) of the satellite SPOT5 scheduling problem is proposed by~\cite{BensanaE1999}. Based on this data, a weighted acyclic digraph model is formulated by~\cite{GabrelV2002}, and solved with a label-setting shortest path algorithm. Alternatively, formulations as generalized knapsack problems can be solved with a tabu search algorithm~\citep{VasquezM2001} or a genetic algorithm~\citep{MansourMAA2010}. Two 0-1 linear programming models are considered by~\cite{GabrelV2006}. Based on the valid inequalities that arise from node packing and the 3-regular independence system polyhedra, a strengthened formulation for the SPOT5 daily photograph scheduling is presented by~\cite{RibeiroGM2010}. However, the benchmark instances are provided without consideration of the constraints that are imposed by a limited observation time of the target. 

\cite{WolfeWJ2000} propose a greedy algorithm and a genetic algorithm based on the assumption that there are only one resource and one observation window for every mission. A single-satellite single-orbit scheduling problem is addressed with a tabu search heuristic in~\citep{CordeauJF2005}, an adaptive meta-heuristic in~\cite{LiuX2017}, and a 0/1 linear programming model in~\citep{SunB2010}. Another 0/1 model based on preprocessing the observation segments is discussed by~\cite{JangJ2013}. The problem of maximizing the total amount of downloaded data is addressed with a mixed-integer programming model and an iterative algorithm~\citep{SpangeloS2015}. There are also several publications that treat the single satellite scheduling as a machine scheduling problem with constraints of operating time windows. The problem is then solved by a heuristic~\citep{BarbulescuL2004,LinW2005,ChengTCE2008,TangpattanakulP2015}. By considering the setup time between two consecutive observations, \cite{LemaitreM2002} introduce the selecting and scheduling problem for an agile Earth observation satellite. \cite{Dilkina2005} take the limited time window and transition time constraints into account. 

In comparison to the single satellite scheduling problem, the use of multiple satellites gives more flexibility and is thus more challenging~\citep{SpangeloS2015}. \cite{WuG2013,ZhangZ2014,XiongJ2016} use graph representations to formulate the problem, for which dynamic programming and ant colony optimization algorithms are proposed to produce a near-optimal solution. To this end, simple sequential missions with conflicts can easily be represented as graphs. However, if the problem involves multiple satellites, the visibility fields of different resources may overlap. Furthermore, several targets may be in the field of view of the same resource simultaneously, and a target may be observed by more than one resource at the same time. Thus, the visible time windows are highly overlapping during the scheduling period, making the combinational characteristic of the problem more prominent. This ultimately renders the uniform modelling of the problem difficult~\citep{YaoF2010}. 

In order to decrease the complexity of the problem and improve computational efficiency, the multi-satellite scheduling problem is often decomposed into the primary problem of mission assignment and the sub-problem of single satellite scheduling~\citep{YaoF2010,WuG2012}. However, since each mission can be observed by multiple resources and since visible time windows interact, the decomposition approach is likely to become trapped in a local optimum of low quality. As a result, a series of mission merging strategies are studied~\cite{XuY2010,WangJ2015}, and the multi-satellite scheduling for dynamic emergency missions is investigated~\cite{NiuX2015,WangJ2015,WangM2014}.
A multitude of different approaches have been investigated to address the problem, such as heuristics~(e.g. greedy approaches~\citep{BianchessiN2008,WangP2011} or local searches~\citep{BonissonePP2006}) and meta-heuristic algorithms~(e.g.\ tabu searches~\citep{HabetD2010,VasquezM2003,BianchessiN2007}, genetic algorithms~\citep{MaoT2012,SunB2010,XhafaF2012,XhafaF2013}, evolutionary algorithms~\citep{BonissonePP2006,SalmanAA2015}, and simulated annealing algorithms~\citep{YaoF2010,XhafaF2013}). 
While these optimization techniques show improvements towards obtaining the optimal or near-optimal solutions, they typically require extensive parameter tuning and cannot provide quality guarantees for the obtained solutions.
\section{A Mixed Integer Linear Programming Model}\label{Sec:Model}
\subsection{Problem Description}
In this section, we describe our mathematical programming model for the multi-satellite scheduling problem. We consider each mission as a point target that has to be observed continuously for a specified time by one of the resources of a satellite. An observation has to be carried out at a certain swing angle and a rotation angle of the resource. Therefore, a setup time between two consecutive observations of the same resource has to be taken into account. A mission requests a certain \emph{imaging type} (visible, multispectral, infrared, or synthetic aperture radar), which must be provided by the corresponding resource. Of course, for modelling the problem, it is not necessary to know the location of the targets. We can determine beforehand at which times a target is visible to each resource and then schedule its observation accordingly. We first define the basic concepts that are used in the following. 

\begin{definition}\label{Def:VTW}
\emph{Visible Time Window}. A time interval during which the target is visible for the resource. 
\end{definition}

\begin{definition}\label{Def:OTW}
\emph{Observation Time Window}. A time interval during which a given resource is assigned to a mission in the scheduling scheme. 

\end{definition}
\begin{definition}\label{Def:FTI}
\emph{Feasible Time Interval}. A continuous time interval that is generated by the union of the overlapping visible time windows of a given resource.
\end{definition}

\begin{definition}\label{Def:CDI}
\emph{Conflict Degree}. The number of candidate missions that can be assigned to the same feasible time interval or subinterval. 
\end{definition}

\noindent
The definitions are described in detail in Appendix~\ref{App:FTI}. Based on these concepts, we introduce the notation used in the following.

\subsection{Notation}

\noindent
As the \emph{scheduling period}, we denote the time interval during which observations can be scheduled, and write $[\SBeg,\SEnd]$, where $\SBeg\ge0$.

\noindent
The set of \emph{missions} is denoted as ${\cal M}=\{M_1,M_2,\ldots,M_n\}$. Each such mission~$M_i$ is specified by its earliest possible observation time~$\ear_i$, its latest possible observation time~$\lat_i$, and the requested duration~$\dur_i$ of the observation. Therefore, to satisfy the mission, a subinterval of length~$\dur_i$ has to be chosen in the time interval $[\ear_i,\lat_i]$. Every mission $M_i$ has a positive \emph{weight}~$\wght_{i}$ measuring its importance (the larger $\wght_i$, the more important the mission).

\noindent
${\cal R}=\{R_1,R_2,\ldots,R_l\}$ is the set of \emph{resources} (cameras, sensors) available on the various satellites. The maximum possible usage time of the resource~$R_j$ in the scheduling period is denoted as~$A_j$. 

\noindent
The availability of resources for missions (in the scheduling period) is specified by appropriate visible time windows. For every resource~$R_j$ and every mission~$M_i$, there is a set $\TW_{ij}=\{\tw_{ij}^1,\ldots,\tw_{ij}^{n_{ij}}\}$ of $n_{ij}$ \emph{visible time windows} during which the resource can be used continuously for the mission. By determining the union of all overlapping visible time windows on resource~$R_j$ over the entire scheduling period, we can compute several disjoint \emph{feasible time intervals} that can be assigned to missions. We denote them as $\RTW_j=\{\rtw_{j}^1,...,\rtw_{j}^{F_j}\}$ where $F_j$ is the number of feasible time intervals of the resource~$R_j$. It is obvious that, for all $\tw_{ij}^k$, $M_i\in M(R_j)$, $k\in{\cal N}_{ij}=\{1,2,...,n_{ij}\}$, we have $\tw_{ij}^k\subseteq \RTW_j$. Each time window $\tw_{ij}^k$ is then given as $\tw_{ij}^k=[\WBeg_{ij}^{k},\WEnd_{ij}^{k}]$.

\noindent
For a mission $M_i\in {\cal M}$, let $R(M_i)\subseteq {\cal R}$ be the set of resources that can be used for this mission. Let $M(R_j)\subseteq {\cal M}$ be the set of missions that a resource $R_j\in {\cal R}$ can service.

\begin{figure}[htbp]
  \centering
  \begin{minipage}{0.38\textwidth}
    \centering
    \includegraphics[height=40mm]{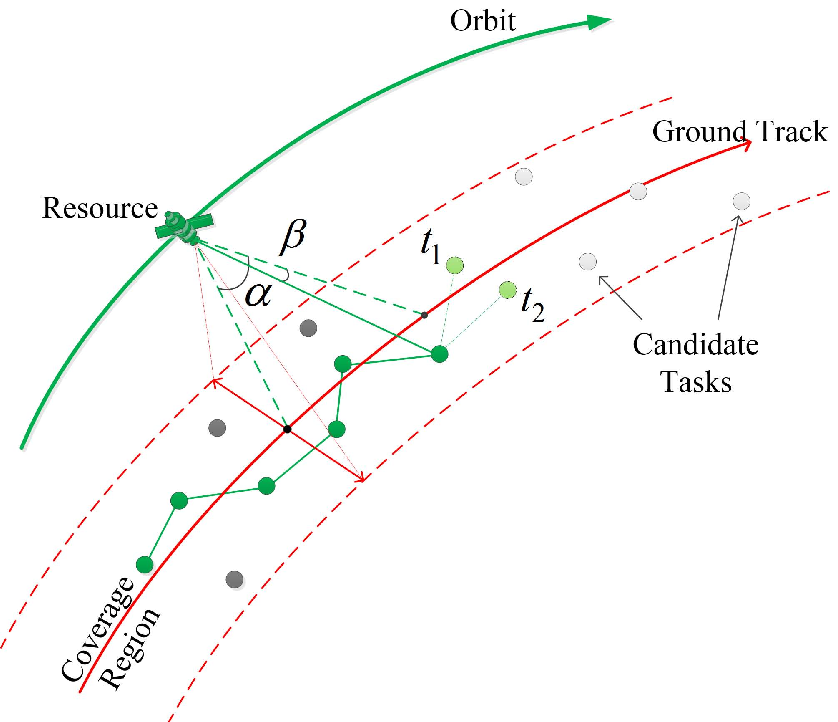}
    \caption{ Resource observation.}
    \label{fig:ResourceObservation}
  \end{minipage}
  \begin{minipage}{0.60\textwidth}
    \centering
    \includegraphics[height=40mm]{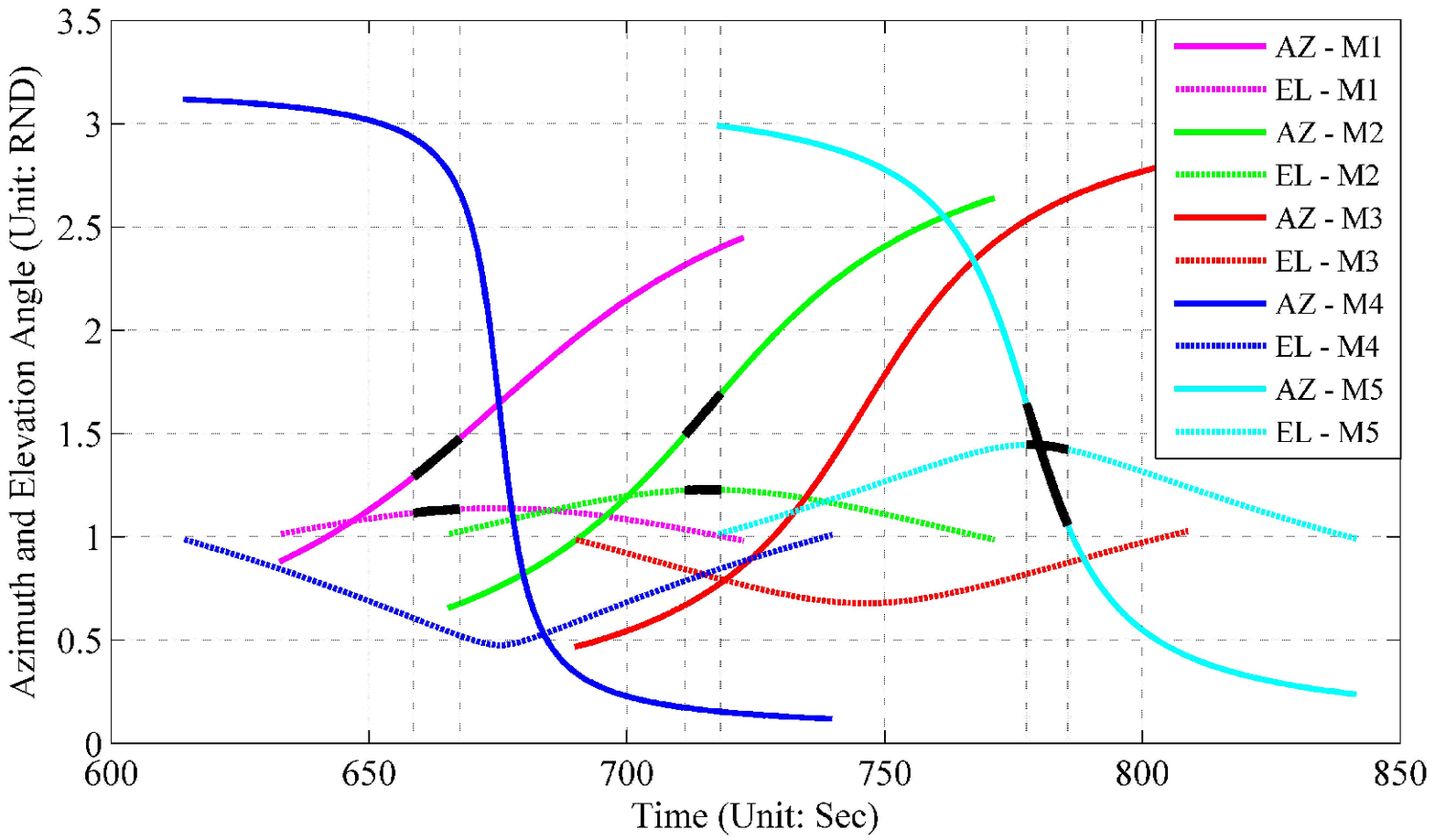}
    \caption{Variation of swing and rotation angles.}
    \label{fig:AzEl}
  \end{minipage}
\end{figure}

\noindent
The resource observation status is illustrated in Fig.~\ref{fig:ResourceObservation}. If resource~$R_j$ is used for observing mission~$M_i$, then this has to happen with a certain \emph{swing angle}~$\alpha$ and a \emph{rotation angle}~$\beta$. These angles are not constant but depend on the position of the resource and thus on the time~$t_i$ at which the observation for $M_i$ starts within one of its feasible time windows. There are existing functions for computing these angles, which we denote as $\alpha_{ijt_i}$ and $\beta_{ijt_i}$, respectively. A schematic overview of the functions is depicted in Fig.~\ref{fig:AzEl}. Here, the elevation angle EL-$M_i$ and azimuth AZ-$M_i$ show angles $\alpha_{ijt_i}$ and $\beta_{ijt_i}$ for a mission $M_i$ and the same resource~$R_j$, depending on the start time $t_i$ of the observation. 

If two consecutive missions are to be carried out by the same resource, then a \emph{setup time} has to be taken into account. The swing angle of resource $R_j$ can be changed by~$\theta_j$ per second and, similarly, its rotation angle by~$\varphi_j$ per second. Furthermore, some time~$\delta_j$ is needed for stabilizing the resource after the angles are adjusted. Angles cannot be adapted simultaneously, so if $M_{i}$ and $M_{i'}$ are two consecutive missions for $R_j$, then the time needed for changing to the correct position is

\begin{equation*}
\mu_{it_ii't_{i'}}^{j}=\frac{|\alpha_{i'jt_{i'}}-\alpha_{ijt_i}|}{\theta_{j}}+\frac{|\beta_{i^{'}jt_{i'}}-\beta_{ijt_i}|}{\varphi_{j}}+\delta_{j}.
\end{equation*}

Of course, $\mu_{it_ii't_{i'}}^{j}\le t_{i'}-t_i-\dur_i$ has to be satisfied.

We simplify the consideration of setup times by computing upper bounds. As we have $\beta_{ij}\in [0,2\pi)$ and $\alpha_{ij}\in [-\alpha_{j}, \alpha_{j}]$ for some maximum angle~$\alpha_j$ depending on the resource, the maximum possible setup time between missions $M_{i}$ and $M_{i'}$ for $R_{j}$ is $\delt_{ii'}^{j}=\frac{2\alpha_{j}}{\theta_{j}}+\frac{\pi}{\varphi_{j}}+\delta_{j}$.


\subsection{Decision Variables}

\noindent
For every mission~$M_i$ and every resource~$R_j$, the binary variable $x_{ij}^{k}$ specifies whether one of the available visible time windows~(${\cal N}_{ij}$) is selected. We let $x_{ij}^{k}=1$ if the visible time window $\tw_{ij}^k$ is used. Otherwise, let $x_{ij}^{k}=0$. 
For every mission~$M_i$, we also have a continuous variable $\OBeg_i$ that denotes the starting time of its observation.


\subsection{Objective}

Assuming that the resources are limited and that not all missions can be carried out, our objective is to schedule either as many missions as possible, i.e.,
\begin{equation}\nonumber\label{model:Objfunction1}
\max \sum_{M_i\in {\cal M}}\;\sum_{R_j\in R(M_i)}\;\sum_{k\in {\cal N}_{ij}}\,x_{ij}^k,
\end{equation}
or to maximize the total weight of accomplished missions, i.e.,
\begin{equation}\nonumber\label{model:Objfunction2}
\max \sum_{M_i\in {\cal M}}\;\sum_{R_j\in R(M_i)}\;\sum_{k\in {\cal N}_{ij}}\wght_i\,x_{ij}^k.
\end{equation}


\subsection{Constraints}

In the following, let $U$ denote a large number depending on the scheduling period (i.e. it serves as the ``Big-M'' required for modelling logical implications).

\paragraph{Mission Accomplishment} It is unlikely that every mission can be carried out. Therefore, although one target can be observed by several resources, the profit of each target counts at most once. For every mission $M_i\in {\cal M}$, we thus have 
\begin{equation}\nonumber\label{model:Constraint_Res1}
\sum_{R_j\in R(M_i)}\;\sum_{k\in{\cal N}_{ij}} x_{ij}^k\;\le\;1.
\end{equation}

\paragraph{Maximum Usage Time} The total observation time in the scheduling period that can be scheduled for a resource is bounded by the given maximum observation time. So for every $R_j\in {\cal R}$, we have the inequality 
\begin{equation}\nonumber\label{model:Constraint_Usg}
\sum_{M_i\in M(R_j)}\;\sum_{k\in{\cal N}_{ij}} \dur_i\cdot x_{ij}^k\;\le\;\usg_j.
\end{equation}

\paragraph{Feasibility of Observation Time} For each $M_i\in{\cal M}$ we have 
\begin{equation}\nonumber\label{model:Constraint_Nonneg}
\begin{split}
\OBeg_i &\;\ge\;\SBeg\\
\OBeg_i+\dur_i &\;\le\;\SEnd.
\end{split}
\end{equation}

\paragraph{Observation Window} If resource $R_j$ and time window $\tw_{ij}^k$ for mission $M_i$ have been selected, then the observation activity has to be placed completely within this interval. This is modeled by the following constraints for all $M_i\in {\cal M}$ and every $R_j\in R(M_i)$, $k\in{\cal N}_{ij}$: 
\begin{equation}\nonumber\label{model:Constraint_Obs}
\begin{split}
\begin{aligned}
\OBeg_i-\WBeg_{ij}^k\cdot x_{ij}^k &\;\ge\;0,\\[0.1cm]
\OBeg_i-(\WEnd_{ij}^k-\dur_i)\cdot x_{ij}^k - U\cdot(1-x_{ij}^k)&\;\le\;0.
\end{aligned}
\end{split}
\end{equation}

\begin{figure}[ht]
  \centering
  \includegraphics[height=45mm]{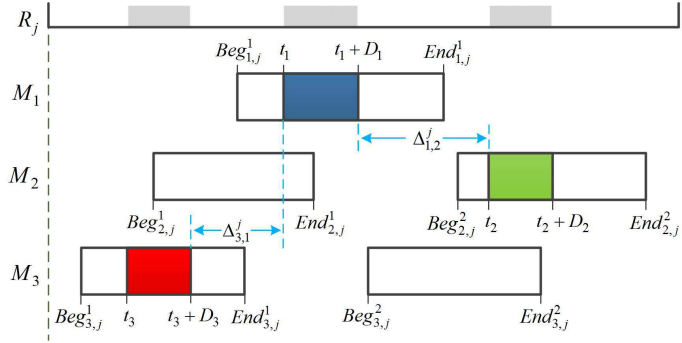}
  \caption{ Observation time window constraint.}
  \label{fig:TimewindowConstraint}
\end{figure}

\paragraph{Setup Time} A minimum transition time for achieving the correct position has to be considered between each pair of consecutive observation activities of the same resource. This situation is depicted in Fig.~\ref{fig:TimewindowConstraint}.

Thus, for all $R_j\in {\cal R}$ and any pair of observations $M_{i},M_{i'}\in {\cal M}(R_j)$, if both missions $M_{i}$ and $M_{i'}$ have been assigned to be carried out by $R_{j}$ then either $\OBeg_{i} \ge \OBeg_{i'} + \dur_{i'} + \delt_{ii'}^{j}$ or $\OBeg_{i'} \ge \OBeg_{i} + \dur_{i} + \delt_{ii'}^{j}$ has to hold.

For modelling this constraint, we introduce binary variables $f_{ii'}^j$ and $f_{i'i}^j$, where $f_{ii'}^j=1$ ($f_{i'i}^j=1$) if and only if both $M_{i}$ and $M_{i'}$ are carried out by $R_j$ and $M_{i}$ is observed after $M_{i'}$ ($M_{i'}$ is observed after~$M_{i}$). Since $f_{ii'}^{j}+f_{i'i}^{j}=\sum_{k}{x_{ij}^{k}}\cdot \sum_{k'}{x_{i'j}^{k'}}$, the disjunction can then be expressed as

\begin{equation}\nonumber\label{model:Constraint_transtime}
\begin{split}
\OBeg_{i} - \OBeg_{i'} &\ge (\dur_{i'} + \delt_{ii'}^{j})\cdot f_{ii'}^j - (U-\dur_{i'})\cdot(1-f_{ii'}^j)\\[0.2cm]
\OBeg_{i'} - \OBeg_{i} &\ge (\dur_{i} + \delt_{ii'}^{j})\cdot f_{i'i}^j - (U-\dur_{i})\cdot(1-f_{i'i}^j).
\end{split}
\end{equation}

For a consistent setting of the new variables, we need 

\begin{equation}\nonumber
\begin{split}
f_{ii'}^j + f_{i'i}^j &\;\le\; \sum_{k\in{\cal N}_{ij}} x_{ij}^{k},\\
f_{ii'}^j + f_{i'i}^j &\;\le\; \sum_{k'\in{\cal N}_{i'j}} x_{i'j}^{k'},\\
f_{ii'}^j + f_{i'i}^j &\;\ge\; \sum_{k\in{\cal N}_{ij}} x_{ij}^{k} + \sum_{k'\in{\cal N}_{i'j}} x_{i'j}^{k'} - 1,\\
\sum_{R_j\in R(M_{i})\cap R(M_{i'})}\;(f_{ii'}^j + f_{i'i}^j)&\;\le\;1.
\end{split}
\end{equation}

\paragraph{Resource Feasibility} 
Motivated by the definitions of the \emph{feasible time interval} and the \emph{conflict degree}~(also see Appendix~\ref{App:FTI}), we consider resource feasibility as a constraint. For a \emph{feasible time interval}~$\rtw_{j}^k$ on~$R_j$, all \emph{effective feasible time subintervals}~$ \srtw_{j}^{kl}\subset\rtw_{j}^k$ are also considered according to the contention conflict degree. In fact, it is a continuous time subinterval within which the visible time windows are highly overlapping and not all of the corresponding candidate missions can be accomplished in $\srtw_{j}^{kl}$. We calculate how many missions can at most be assigned to a subinterval~$\srtw_{j}^{kl}$, and denote them as $\srn_{j}^{kl}$. For each such effective feasible time subinterval~$\srtw_{j}^{kl}$, we have inequalities 

\begin{equation}\nonumber\label{model:Constraint_Res2}
\begin{split}
\sum_{M_i\in {\cal M}(R_j)}\;\sum_{k'\in{\cal N}_{ij}} x_{ij}^{k'}\;& \le\;\srn_{j}^{kl} \\
\sum_{M_i\in {\cal M}(R_j)}\;\sum_{k'\in{\cal N}_{ij}} (\dur_i+\delta_j)\cdot x_{ij}^{k'}\;& \le\;|\srtw_{j}^{kl}|+\delta_j ,
\end{split}
\end{equation}
where $tw_{ij}^{k'}\subset\srtw_{j}^{kl}$.

The generation of the effective feasible time subinterval~$\srtw_{j}^{kl}$ and the computation of the corresponding maximum assignment capacity~$\srn_{j}^{kl}$ are described in detail in Appendix~\ref{App:SRTW}. All operations with respect to the resource feasibility are computed in preprocessing, and the computation results are regarded as constraints in modelling. Furthermore, these newly proposed formulations are a set of effective inequalities and produce a significant improvement in obtaining a tighter upper bound of instances. 

\paragraph{Integrality Constraints} For all $M_i\in{\cal M}$, all $R_j\in {\cal R}(M_i)$, and all $k\in{\cal N}_{ij}$ we have 
\begin{equation}\nonumber\label{model:Constraint_Int}
x_{ij}^k\;\in\;\{0,1\}. 
\end{equation}

\subsection{Improved Constraints}
\begin{figure}[!htbp]
  \centering
  \includegraphics[width=99mm]{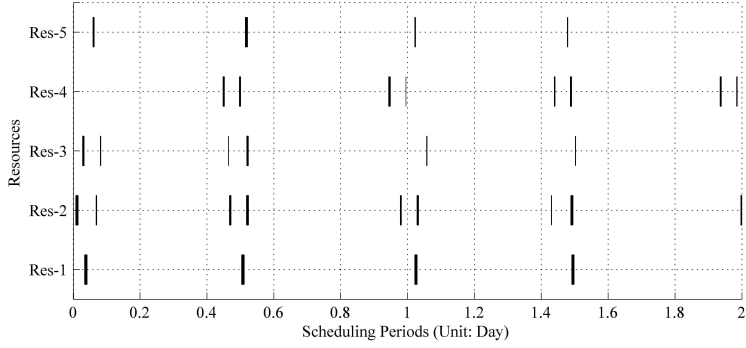}
  \caption{ Distribution of feasible time intervals.}
  \label{fig:FeasibleTimeInterval4}
\end{figure}

A mission may have several disjoint visible time windows for each resource and the fraction of the time window that is needed for observation may be comparatively small. Fig.~\ref{fig:FeasibleTimeInterval4} illustrates the distribution of all feasible time intervals over time. We find that the availability distribution of resources is quite sparse over the entire scheduling period. 

In the first model, the value of $\OBeg_i$ is always feasible for the entire scheduling period due to the introduction of~$U$ as ``Big-M''. Typically, for large values of $U$, this is an obstacle for obtaining a good bound from the LP relaxation. Thus, we avoid the use of~$U$ by defining new (continuous) variables~$\OBeg_{ij}^{k}$ that denote the starting time of mission~$M_i$ if $x_{ij}^{k}$ is selected. As bounds we then have $\OBeg_{ij}^k\ge\SBeg$ and $\OBeg_{ij}^k+\dur_i\le\SEnd$. Based on this deliberation, new constraints are introduced as follows.

\paragraph{Observation Window} If the resource $R_j$ and a corresponding visible time window $\tw_{ij}^k$ for mission $M_i$ have been selected, then the observation has to be completely placed within this interval. For all $M_i\in {\cal M}$ and $R_j\in {\cal R}(M_i)$, $k\in{\cal N}_{ij}$, we thus have 

\begin{equation}\nonumber\label{model:Constraint_ObsE}
\begin{split}
\begin{aligned}
\OBeg_{ij}^k-\WBeg_{ij}^k\cdot x_{ij}^k &\;\ge\;0,\\[0.1cm]
\OBeg_{ij}^k-(\WEnd_{ij}^k-\dur_i)\cdot x_{ij}^k &\;\le\;0.
\end{aligned}
\end{split}
\end{equation}

\paragraph{Setup Time} For all resources $R_j$ and any pair of observations $M_{i},M_{i'}\in {\cal M}(R_j)$, if both missions $M_{i}$ and $M_{i'}$ have been assigned to be carried out by $R_{j}$, then we obtain tighter constraints for the respective starting times since either
$(\dur_{i'} + \delt_{ii'}^{j})\le \OBeg_{ij}^{k} - \OBeg_{i'j}^{k'}\le \WEnd_{ij}^{k}-\dur_{i}-\WBeg_{i'j}^{k'}$
or
$(\dur_{i} + \delt_{ii'}^{j})\le \OBeg_{i'j}^{k'} - \OBeg_{ij}^{k}\le \WEnd_{i'j}^{k'}-\dur_{i'}-\WBeg_{ij}^{k}$.
Here, we distinguish three cases.

\begin{enumerate}
\item
For resource $R_j$, if the candidate visible time windows $\tw_{ij}^{k}$ and $\tw_{i'j}^{k'}$ are potentially ``ordered'' and satisfy that $\WEnd_{ij}^{k}-\dur_{i}-\WBeg_{i'j}^{k'}<\delt_{ii'}^{j}$ (that is, if both of them are assigned, then mission~$M_{i}$ will definitely be executed before mission~$M_{i'}$), we have

\begin{equation}\nonumber\label{model:Constraint_transtimeE1}
\OBeg_{i'j}^{k'} - \OBeg_{ij}^{k} \ge (\dur_{i} + \delt_{ii'}^{j})\cdot x_{ij}^{k} - (\WEnd_{ij}^{k} + \delt_{ii'}^{j})\cdot(1 - x_{i'j}^{k'}).
\end{equation}

\item
For resource $R_j$, if the candidate visible time windows $\tw_{i'j}^{k'}$ and $\tw_{ij}^{k}$ are potentially ``ordered'' and satisfy that $\WEnd_{i'j}^{k'}-\dur_{i'}-\WBeg_{ij}^{k}<\delt_{ii'}^{j}$ (that is, if both of them are assigned, then mission~$M_{i'}$ will definitely be executed before mission~$M_{i}$), we have

\begin{equation}\nonumber\label{model:Constraint_transtimeE2}
\OBeg_{ij}^{k} - \OBeg_{i'j}^{k'} \ge (\dur_{i'} + \delt_{ii'}^{j})\cdot x_{i'j}^{k'} - (\WEnd_{i'j}^{k'} + \delt_{ii'}^{j})\cdot(1 - x_{ij}^{k}).
\end{equation}

\item
Only if the candidate visible time windows $\tw_{ij}^{k}$ and $\tw_{i'j}^{k'}$ are overlapping we introduce binary variables $f_{jii'}^{kk'}$ and $f_{ji'i}^{k'k}$, where $f_{jii'}^{kk'}=1$ if $M_{i}$ is observed after $M_{i'}$ and $f_{ji'i}^{kk'}=1$ if $M_{i'}$ is observed after~$M_{i}$ (by $R_j$). Since $f_{jii'}^{kk'}+f_{ji'i}^{k'k}=x_{ij}^{k}\cdot x_{i'j}^{k'}$, the disjunction can be expressed as

\begin{equation}\nonumber\label{model:Constraint_transtimeE3-1}
\begin{split}
\OBeg_{ij}^{k} - \OBeg_{i'j}^{k'} &\ge (\dur_{i'} + \delt_{ii'}^{j})\cdot f_{jii'}^{kk} - (\WEnd_{i'j}^{k'}-\dur_{i'}-\WBeg_{ij}^{k}\cdot f_{ji'i}^{k'k})\cdot(1 - f_{jii'}^{kk'}),\\[0.2cm]
\OBeg_{i'j}^{k'} - \OBeg_{ij}^{k} &\ge (\dur_{i} + \delt_{ii'}^{j})\cdot f_{ji'i}^{k'k} - (\WEnd_{'j}^{k}-\dur_{i}-\WBeg_{i'j}^{k'}\cdot f_{jii'}^{kk'})\cdot(1 - f_{ji'i}^{k'k}),
\end{split}
\end{equation}

or equivalently as

\begin{equation}\nonumber\label{model:Constraint_transtimeE3-2}
\begin{split}
\OBeg_{ij}^{k} - \OBeg_{i'j}^{k'} &\ge (\WEnd_{i'j}^{k'} + \delt_{ii'}^{j})\cdot f_{jii'}^{kk'} + \WBeg_{ij}^{k}\cdot f_{ji'i}^{k'k} - (\WEnd_{i'j}^{k'}-\dur_{i'}),\\[0.2cm]
\OBeg_{i'j}^{k'} - \OBeg_{ij}^{k} &\ge (\WEnd_{ij}^{k} + \delt_{ii'}^{j})\cdot f_{ji'i}^{k'k} + \WBeg_{i'j}^{k'}\cdot f_{jii'}^{kk'} - (\WEnd_{ij}^{k}-\dur_{i}).
\end{split}
\end{equation}

For the consistent setting of the new binary variables we need

\begin{equation}\nonumber
\begin{split}
f_{jii'}^{kk'} + f_{ji'i}^{k'k} &\le x_{ij}^{k}\\
f_{jii'}^{kk'} + f_{ji'i}^{k'k} &\le x_{i'j}^{k'}\\
f_{jii'}^{kk'} + f_{ji'i}^{k'k} &\ge x_{ij}^{k} + x_{i'j}^{k'} - 1\\
\sum_{R_j\in {\cal R}(M_{i})\cap {\cal R}(M_{i'})}\;\sum_{k\in{\cal N}_{ij}}\;\sum_{k'\in{\cal N}_{i'j}}\;(f_{jii'}^{kk'} + 
f_{ji'i}^{k'k})&\le 1.
\end{split}
\end{equation}

With the introduction of these variables, we do not have~$U$~(serving as ``Big-M'') anymore. Instead, additional 5-index variables are introduced for pairs of overlapping visible time windows. However, for real-world problem instances where there are many visible time windows for the missions, the number of additional variables is acceptable. We give a detailed analysis in subsection~\ref{Subsec:ModelInfo}. The newly introduced constraints also help to reformulate the conflict segments pairs, hence decreasing the number of binary variables in the model~\citep{JangJ2013}, while simultaneously eliminating the ``Big-M'' in formulating the satellite range scheduling problem~\citep{LuoK2017}. 
\end{enumerate}

\section{Computational Experiments}\label{Sec:Simulation}

In the following, we describe our experiments on several test instances.

\subsection{Test Instances}

For analyzing the performance of our model, we generate several test instances. We use the current on-orbit environment and disaster monitoring satellites HJ-1A, HJ-1B, and HJ-1C, which can carry out large-scale, all-weather, and 24h dynamic monitoring for the ecological environment and disaster. The satellite HJ-1A is equipped with a CCD scanner and a hyperspectral imager, HJ-1B is equipped with a CCD scanner and an infrared scanner, and HJ-1C is equipped with an s-wave band synthetic aperture radar that has two working modes. The satellites HJ-1A and HJ-1B are located about 650km above the Earth's surface in a sun-synchronous orbit with 14.737 orbits per day. The satellite HJ-1C is located about 500km above the Earth's surface in a sun-synchronous orbit with 15.22 orbits per day. The resources on satellites HJ-1A and HJ-1B can observe the complete surface of the Earth in 2 days. Three types R, M, and C of spot target sets are generated with different resource conflict characteristics.

\begin{itemize}
\item[R] Targets are generated randomly and uniformly distributed over the entire land-area on Earth.
\item[C] All targets are randomly generated in clusters over the entire land-area on Earth.
\item[M] Targets with a high number of conflicts are generated manually and clustered over several regions on Earth.
\end{itemize}

The observation times~$\dur_i$ for the missions are integers that are generated uniformly from the interval $[3,10]$. If weights~$\wght_{i}$ are taken into account, they are integers that are generated uniformly from $[1,10]$.

Table~\ref{tab:TestInstance} shows a summary of the generated test instances. The hypothetical start time of the scheduling period is 2016-06-01 06:00:00. The scheduling horizon is 24 or 48~hours. By combining different satellites and target sets, 37 problem instances are generated. 

\begin{table}[!htbp]
\scriptsize 
\caption{Test Instances}
\label{tab:TestInstance}
\centering
\linespread{1.0}\selectfont
\begin{tabular}{rrccr|rrccr}
\hline\noalign{\smallskip}
Instance &
Period & 
$|\cal R|$ & 
$|\cal M|$ & 
$\sum w_i$ & 
Instance &
Period & 
$|\cal R|$ & 
$|\cal M|$ & 
$\sum w_i$ \\
\noalign{\smallskip}\hline\noalign{\smallskip}
M-1 & 24H & 3 & C100 & 621 & C-1 & 24H & 3 & M100 & 620 \\
M-2 & 24H & 3 & C200 & 1191 & C-2 & 24H & 3 & M200 & 1218 \\
M-3 & 24H & 3 & C300 & 1790 & C-3 & 24H & 3 & M300 & 1785 \\
M-4 & 24H & 5 & C100 & 621 & C-4 & 24H & 5 & M100 & 620 \\
M-5 & 24H & 5 & C200 & 1191 & C-5 & 24H & 5 & M200 & 1218 \\
M-6 & 24H & 5 & C300 & 1790 & C-6 & 24H & 5 & M300 & 1785 \\
M-7 & 24H & 5 & C400 & 2401 & C-7 & 24H & 5 & M400 & 2385 \\
M-8 & 24H & 5 & C500 & 2976 & C-8 & 24H & 5 & M500 & 2998 \\
M-9 & 48H & 5 & C100 & 621 & C-9 & 48H & 5 & M100 & 620 \\
M-10 & 48H & 5 & C200 & 1191 & C-10 & 48H & 5 & M200 & 1218 \\
M-11 & 48H & 5 & C300 & 1790 & C-11 & 48H & 5 & M300 & 1785 \\
M-12 & 48H & 5 & C400 & 2401 & C-12 & 48H & 5 & M400 & 2385 \\
M-13 & 48H & 5 & C500 & 2976 & C-13 & 48H & 5 & M500 & 2998 \\
R-1 & 24H & 4 & R300 & 1821 & R-7 & 24H & 6 & R600 & 3617 \\
R-2 & 24H & 4 & R400 & 2427 & R-8 & 24H & 6 & R700 & 4193 \\
R-3 & 24H & 4 & R500 & 3004 & R-9 & 24H & 6 & R800 & 4783 \\
R-4 & 24H & 6 & R300 & 1821 & R-10 & 24H & 6 & R900 & 5371 \\
R-5 & 24H & 6 & R400 & 2427 & R-11 & 24H & 6 & R1000 & 5987 \\
R-6 & 24H & 6 & R500 & 3004 &  &  &  &  & \\
\hline
\end{tabular}
\end{table}

Further details on all instances are provided in Appendix~\ref{App:InstAnaly}. By computing some characteristic numbers, we quantify the availability of resources and the complexity of instances. We calculate the maximum possible number of missions~($\rn$) that can be assigned to a resource according to constraint~\eqref{model:Constraint_Res2}. It is given as the sum of $\rn_j^k$ of feasible time intervals of the resource~$R_j$. Furthermore, we introduce the resource contention degree~($\conf$), which reflects the limitation of each feasible time interval of resources, as well as the average conflict degree of missions, which is calculated as $\conf_j=\frac{T_j-F_j}{F_j}$. Intuitively, it denotes how many missions can be assigned synchronously to a resource during the entire scheduling period, while $\conf_j=0$ indicates that there is no resource contention conflict. These two conflict indicators illustrate the potential complexity of the instance when assigning resources to missions.

In comparison to $N$, $\rn$ decreases quickly for higher values of $\conf$, indicating a more exact upper bound. Additionally, the potential assignment opportunity of missions denotes the flexibility in assigning a resource along with an observation time to a mission (including the average number of visible time windows of a mission~($\paon$) and the average visible time of a mission~($\paot$)). Here, it is obvious that even if two missions have the same total visible time duration, it is more difficult to assign the mission that has the higher number of visible time windows. 

\subsection{Comparison of Proposed Models}\label{Subsec:ModelInfo}
A summary of the decision variables used in the different formulations is shown in Table~\ref{tab:DecisionVariableInfo}.

\begin{table}[ht]
\scriptsize 
\caption{Model Infomation}
\label{tab:DecisionVariableInfo}
\centering
\linespread{1.5}\selectfont
\begin{tabular}{|c|c|l|l|}
\hline
\multirow{4}{1.0cm}{MILP} & \multirow{2}{1.2cm}{Decision Variables} & $x_{ij}^k\,\in\,\{0,1\}$ & $M_i\in{\cal M}$, all $R_j\in {\cal R}(M_i)$, and all $k\in{\cal N}_{ij}$ \\
\cline{3-4}
 &  & $t_i\,\in\,\mathbb{R}$ & $M_i\in{\cal M}$ \\
\cline{2-4}
 & \multirow{2}{1.2cm}{Introduced Variables} & $f_{ii'}^j\,\in\,\{0,1\}$ & \multirow{2}{*}{$M_i,M_{i'}\in{\cal M}$, all $R_j\in {\cal R}(M_{i})\cap {\cal R}(M_{i'})$} \\
\cline{3-3}
 &  & $f_{i'i}^j\,\in\,\{0,1\}$ &  \\
\hline
\multirow{5}{1.0cm}{Improved MILP} & \multirow{2}{1.2cm}{Decision Variables} & $x_{ij}^k\,\in\,\{0,1\}$ & \multirow{2}{*}{$M_i\in{\cal M}$, all $R_j\in {\cal R}(M_i)$, and all $k\in{\cal N}_{ij}$} \\
\cline{3-3}
 &  & $t_{ij}^k\,\in\,\mathbb{R}$ &  \\
\cline{2-4}
 & \multirow{3}{1.2cm}{Introduced Variables} & $f_{jii'}^{kk'}\,\in\,\{0,1\}$ & $M_i,M_{i'}\in{\cal M}$, all $R_j\in {\cal R}(M_{i})\cap {\cal R}(M_{i'})$  \\
\cline{3-3}
 &  & \multirow{2}{*}{$f_{ji'i}^{k'k}\,\in\,\{0,1\}$} & all $k\in{\cal N}_{ij}$ and $k'\in{\cal N}_{i'j}$ \\
 &  &  & $tw_{ij}^{k}$ and $tw_{i'j}^{k'}$ are overlapped. \\
\hline
\end{tabular}
\end{table}

Furthermore, we investigate the effect of the preprocessing step~(i.e., the reduction of the size of the search space) and compare the proposed models over different instances~(i.e., the number of variables and constraints of a model). The results are shown in detail in Table~\ref{tab:SchedulingModelInfo} and indicate an improvement of each phase in the overall performance of the proposed method. 

\begin{table}[ht]
\scriptsize 
\caption{Model Infomation}
\label{tab:SchedulingModelInfo}
\centering
\linespread{1.0}\selectfont
\begin{tabular}{lrrrrrrr}
\hline\noalign{\smallskip}
\multirow{3}{*}{Ins.} & 
\multirow{3}{0.8cm}{\raggedleft $n'$} & 
\multicolumn{3}{c}{\centering MILP} &
\multicolumn{3}{c}{\centering Improved MILP} \\
\cmidrule(l){3-5}\cmidrule(l){6-8}\noalign{\smallskip}
 & & 
 \multicolumn{1}{p{1.1cm}}{\raggedleft $\mVC$} & 
 \multicolumn{1}{p{1.1cm}}{\raggedleft $\mVB$} & 
 \multicolumn{1}{p{1.1cm}}{\raggedleft $\mC$} & 
 \multicolumn{1}{p{1.1cm}}{\raggedleft $\mVC$} & 
 \multicolumn{1}{p{1.1cm}}{\raggedleft $\mVB$} & 
 \multicolumn{1}{p{1.1cm}}{\raggedleft $\mC$} \\
\noalign{\smallskip}\hline\noalign{\smallskip}
\textbf{C-1} & 0 & 100  & 4082  & 10463  & 100 & 1838 & 7713 \\
\textbf{C-2} & 6 & 194  & 15234  & 38620  & 218 & 3456 & 20510 \\
\textbf{C-3} & 2 & 298  & 36466  & 92019  & 360 & 8794 & 53973 \\
\textbf{C-4} & 2 & 98  & 8071  & 20595  & 189 & 3543 & 14114 \\
\textbf{C-5} & 6 & 194  & 38744  & 97831  & 454 & 7176 & 44853 \\
\textbf{C-6} & 2 & 298  & 78881  & 198598  & 701 & 17349 & 106032 \\
\textbf{C-7} & 15 & 385  & 129662  & 326594  & 992 & 27780 & 192205 \\
\textbf{C-8} & 13 & 487  & 198407  & 499268  & 1265 & 45169 & 311769 \\
\textbf{C-9} & 3 & 97  & 9378  & 24157  & 354 & 7280 & 35404 \\
\textbf{C-10} & 7 & 193  & 50217  & 126800  & 667 & 11945 & 92045 \\
\textbf{C-11} & 3 & 297  & 116704  & 293764  & 1098 & 30796 & 235087 \\
\textbf{C-12} & 18 & 382  & 218518  & 550073  & 1664 & 43784 & 463025 \\
\textbf{C-13} & 16 & 484  & 345330  & 868124  & 2086 & 67448 & 724176 \\
\noalign{\smallskip}
\textbf{M-1} & 17 & 83  & 2969  & 7626  & 83 & 1453 & 5964 \\
\textbf{M-2} & 18 & 182  & 13938  & 35388  & 206 & 3308 & 19270 \\
\textbf{M-3} & 27 & 273  & 32290  & 81562  & 320 & 8278 & 47707 \\
\textbf{M-4} & 19 & 81  & 5922  & 15125  & 158 & 2916 & 11260 \\
\textbf{M-5} & 19 & 181  & 31354  & 79275  & 414 & 6738 & 39564 \\
\textbf{M-6} & 29 & 271  & 72947  & 183759  & 643 & 16907 & 97369 \\
\textbf{M-7} & 50 & 350  & 152005  & 382102  & 797 & 16431 & 136066 \\
\textbf{M-8} & 93  & 407  & 182919  & 459547  & 885 & 17681 & 158157 \\
\textbf{M-9} & 19  & 81  & 7977  & 20519  & 313 & 6263 & 30347 \\
\textbf{M-10} & 22  & 178  & 51349  & 129772  & 721 & 13201 & 103415 \\
\textbf{M-11} & 31  & 269  & 123785  & 311638  & 1119 & 31573 & 247473 \\
\textbf{M-12} & 64  & 336  & 221047  & 556076  & 1323 & 26747 & 328150 \\
\textbf{M-13} & 118  & 382  & 283211  & 711799  & 1525 & 29387 & 415919 \\
\noalign{\smallskip}
\textbf{R-1} & 125 & 175  & 13765  & 34892  & 227 & 3321 & 21800 \\
\textbf{R-2} & 181 & 219  & 21159  & 53458  & 273 & 4419 & 30080 \\
\textbf{R-3} & 237 & 263  & 29208  & 73625  & 324 & 4638 & 36312 \\
\textbf{R-4} & 125 & 175  & 28598  & 72309  & 438 & 6206 & 40462 \\
\textbf{R-5} & 181 & 219  & 45265  & 114097  & 533 & 8493 & 57609 \\
\textbf{R-6} & 237 & 263  & 63593  & 159954  & 631 & 8895 & 70923 \\
\textbf{R-7} & 284 & 316  & 92423  & 232641  & 829 & 11651 & 112816 \\
\textbf{R-8} & 305 & 395  & 148393  & 373188  & 1079 & 15525 & 180564 \\
\textbf{R-9} & 347 & 453  & 215280  & 540854  & 1342 & 16974 & 263464 \\
\textbf{R-10} & 379 & 521  & 279252  & 701080  & 1514 & 17806 & 324658 \\
\textbf{R-11} & 389 & 611  & 405523  & 1017699  & 1819 & 22687 & 461934 \\
\noalign{\smallskip}\hline
\end{tabular}
\end{table}

In Table~\ref{tab:SchedulingModelInfo}, $n'$ denotes the number of missions that are scheduled during preprocessing. It is calculated based on the \emph{effective feasible time subinterval}~(see Appendix~\ref{App:SRTW}). With $\mVC$, $\mVB$, and $\mC$, we denote the number of continuous variables, the number of binary variables, and the number of constraints, respectively. The results show that the preprocessing is especially effective for instances with randomly distributed targets~${\cal R}$. In the MILP, the number of continuous variables equals the number of missions, and the number of binary variables exponentially increases with the number of missions, meaning that $\mVC=n-n'$ and $\mVB\approx (n-n')^2$. In contrast, for the improved MILP, the number of continuous variables equals the number of visible time windows, meaning that $\mVC=(n-n')\cdot \paon$. Due to the linearization of the formulation, the binary variable is introduced only if the two candidate visible time windows are overlapping. Thus, compared to the MILP, the number of binary variables is smaller, especially for the instance $\cal C$ of target sets~(see Table~\ref{tab:ResourceUtilisation}), where $\mVB \ll (n-n')^2$. 

\subsection{Optimization Results}

We test our model on all instances with Gurobi 6.5.1 on a 3.40GHz PC with 16GB RAM and 8~cores. The maximum run time is set to 6~hours. Both objective functions, i.e., maximizing the number of scheduled missions as well as maximizing the total weights of scheduled missions, are considered. The results are shown in Table~\ref{tab:SchedulingResultsNumObj} and Table~\ref{tab:SchedulingResultsWegObj}, respectively. For each instance, we show the \emph{root upper bound}, the \emph{final upper bound} found by the Gurobi, and the value of the \emph{best solution found}. If the running time is not exceeded, then the latter two are equal and represent the optimum value. The \emph{gap} is computed as $(\emph{final upper bound}-\emph{best solution})/\emph{final upper bound}$. The CPU time is shown if the optimum is found in less than 6~hours. 

\begin{table}[!ht]
\scriptsize 
\caption{Optimization results: number of missions}
\label{tab:SchedulingResultsNumObj}
\centering
\linespread{1.1}\selectfont
\begin{tabular}{lrrrrrrrr}
\hline\noalign{\smallskip}
\multirow{4}{*}{Ins.} & 
\multicolumn{8}{c}{\centering Objective: Maximize the total number of assigned missions} \\
\cmidrule(l){2-9}\noalign{\smallskip}
 & 
\multicolumn{4}{c}{\centering Produced by the MILP} &
\multicolumn{4}{c}{\centering Produced by the Improved MILP} \\
\cmidrule(l){2-5}\cmidrule(l){6-9}\noalign{\smallskip}
 & 
 \multicolumn{1}{p{1.7cm}}{\raggedleft Root/Final upper bound} & 
 \multicolumn{1}{p{0.8cm}}{\raggedleft Best Result} & 
 \multicolumn{1}{p{0.7cm}}{\raggedleft \multirow{2}{*}{\centering Gap}} & 
 \multicolumn{1}{p{0.6cm}}{\raggedleft Run time(s)} & 
 \multicolumn{1}{p{1.8cm}}{\raggedleft Root/Final upper bound} & 
 \multicolumn{1}{p{0.8cm}}{\raggedleft Best Result} & 
 \multicolumn{1}{p{0.7cm}}{\raggedleft \multirow{2}{*}{\centering Gap}} & 
 \multicolumn{1}{p{0.6cm}}{\raggedleft Run time(s)} \\
\noalign{\smallskip}\hline\noalign{\smallskip}
\textbf{C-1} & 26.00/26 & 26 & 0.00\% & 0 & 26.00/26 & 26 & 0.00\% & 0 \\
\textbf{C-2} & 82.00/82 & 82 & 0.00\% & 19 & 82.00/82 & 82 & 0.00\% & 0 \\
\textbf{C-3} & 85.89/85 & 85 & 0.00\% & 265 & 85.00/85 & 85 & 0.00\% & 1 \\
\textbf{C-4} & 44.00/44 & 44 & 0.00\% & 23 & 44.00/44 & 44 & 0.00\% & 2 \\
\textbf{C-5} & 129.84/128 & 126 & 1.56\% & -- & 130.84/129 & 128 & 0.78\% & -- \\
\textbf{C-6} & 133.00/133 & 125 & 6.02\% & -- & 132.00/131 & 131 & 0.00\% & 683 \\
\textbf{C-7} & 163.98/163 & 157 & 3.68\% & -- & 164.98/163 & 162 & 0.61\% & -- \\
\textbf{C-8} & 181.55/179 & 125 & 30.17\% & -- & 183.85/179 & 178 & 0.56\% & -- \\
\textbf{C-9} & 73.76/72 & 71 & 1.39\% & -- & 73.76/72 & 72 & 0.00\% & 157 \\
\textbf{C-10} & 155.84/155 & 121 & 21.94\% & -- & 157.84/155 & 155 & 0.00\% & 377 \\
\textbf{C-11} & 178.00/178 & 132 & 25.84\% & -- & 178.00/176 & 175 & 0.57\% & -- \\
\textbf{C-12} & 265.79/265 & 160 & \textbf{39.62\%} & -- & 266.79/265 & 259 & 2.26\% & -- \\
\textbf{C-13} & 291.97/290 & 180 & 37.93\% & -- & 293.97/290 & 285 & 1.72\% & -- \\
\noalign{\smallskip}
\textbf{M-1} & 39.00/39 & 39 & 0.00\% & 2 & 39.00/39 & 39 & 0.00\% & 0 \\
\textbf{M-2} & 86.00/86 & 86 & 0.00\% & 15 & 86.00/86 & 86 & 0.00\% & 1 \\
\textbf{M-3} & 108.00/107 & 107 & 0.00\% & 153 & 108.00/107 & 107 & 0.00\% & 8 \\
\textbf{M-4} & 53.00/53 & 53 & 0.00\% & 17 & 53.00/53 & 53 & 0.00\% & 1 \\
\textbf{M-5} & 124.00/124 & 121 & 2.42\% & -- & 123.00/122 & 122 & 0.00\% & 139 \\
\textbf{M-6} & 151.00/151 & 143 & 5.30\% & -- & 150.00/150 & 149 & 0.67\% & -- \\
\textbf{M-7} & 228.93/225 & 197 & 12.44\% & -- & 227.93 /224 & 223 & 0.45\% & -- \\
\textbf{M-8} & 309.81/305 & 282 & 7.54\% & -- & 308.81 /304 & 303 & 0.33\% & -- \\
\textbf{M-9} & 82.92/81 & 81 & 0.00\% & 1581 & 82.92/81 & 81 & 0.00\% & 25 \\
\textbf{M-10} & 162.99/162 & 138 & 14.81\% & -- & 162.99/162 & 161 & 0.62\% & -- \\
\textbf{M-11} & 206.00/206 & 158 & 23.30\% & -- & 206.00/205 & 204 & 0.49\% & -- \\
\textbf{M-12} & 305.88/301 & 218 & 27.57\% & -- & 304.88/301 & 296 & 1.66\% & -- \\
\textbf{M-13} & 402.00/402 & 294 & 26.87\% & -- & 403.99/401 & 391 & \textbf{2.49\%} & -- \\
\noalign{\smallskip}
\textbf{R-1} & 210.92/209 & 209 & 0.00\% & 118 & 209.00/209 & 209 & 0.00\% & 1 \\
\textbf{R-2} & 287.92/286 & 286 & 0.00\% & 113 & 286.00/286 & 286 & 0.00\% & 1 \\
\textbf{R-3} & 381.00/380 & 380 & 0.00\% & 223 & 380.00/380 & 380 & 0.00\% & 1 \\
\textbf{R-4} & 239.82/238 & 236 & 0.84\% & -- & 238.75/237 & 236 & 0.42\% & -- \\
\textbf{R-5} & 317.82/316 & 314 & 0.63\% & -- & 315.94/315 & 314 & 0.32\% & -- \\
\textbf{R-6} & 412.00/411 & 409 & 0.49\% & -- & 410.94/410 & 409 & 0.24\% & -- \\
\textbf{R-7} & 507.00/507 & 434 & 14.40\% & -- & 506.99 /505 & 502 & 0.59\% & -- \\
\textbf{R-8} & 581.83/581 & 497 & 14.46\% & -- & 583.83 /580 & 579 & 0.17\% & -- \\
\textbf{R-9} & 682.83/682 & 593 & 13.05\% & -- & 682.83 /679 & 677 & 0.29\% & -- \\
\textbf{R-10} & 763.83/763 & 670 & 12.19\% & -- & 765.83 /762 & 760 & 0.26\% & -- \\
\textbf{R-11} & --/-- & -- & -- & -- & 839.33 /835 & 833 & 0.24\% & -- \\
\noalign{\smallskip}\hline
\end{tabular}
\end{table}

In Table~\ref{tab:SchedulingResultsNumObj}, we see that in comparison to the results produced by the MILP, a larger number of optimal solutions can be obtained by the improved MILP. For most of the problem instances, the tightest upper bounds can be efficiently generated by the improved MILP. Optimal solutions are usually obtained in less than 1,000~s. Within 6~hours, solutions with a small optimality gap can be determined. The worst gap among all results is~2.49\% for the improved model, whereas the worst generated gap among results produced by the MILP is~39.62\%. The MILP also fails to obtain a feasible solution for instance R-11. In combination with Table~\ref{tab:SchedulingModelInfo}, these findings indicate that the performance of the proposed models varies with the size of the instance. The advantages are similar to those obtained when maximizing the total weights of assigned missions shown in Table~\ref{tab:SchedulingResultsWegObj}. Here, the worst gap produced by the improved MILP is 2.95\%, while the worst gap produced by the MILP is 35.97\%.

\begin{table}[!ht]
\scriptsize 
\caption{Optimization results: weight of missions}
\label{tab:SchedulingResultsWegObj}
\centering
\linespread{1.1}\selectfont
\begin{tabular}{lrrrrrrrr}
\hline\noalign{\smallskip}
\multirow{4}{*}{Ins.} & 
\multicolumn{8}{c}{\centering Objective: Maximize the total weight of assigned missions} \\
\cmidrule(l){2-9}\noalign{\smallskip}
 & 
\multicolumn{4}{c}{\centering Produced by the MILP} &
\multicolumn{4}{c}{\centering Produced by the Improved MILP} \\
\cmidrule(l){2-5}\cmidrule(l){6-9}\noalign{\smallskip}
 & 
 \multicolumn{1}{p{1.7cm}}{\raggedleft Root/Final upper bound} & 
 \multicolumn{1}{p{0.8cm}}{\raggedleft Best Result} & 
 \multicolumn{1}{p{0.7cm}}{\raggedleft \multirow{2}{*}{\centering Gap}} & 
 \multicolumn{1}{p{0.6cm}}{\raggedleft Run time(s)} & 
 \multicolumn{1}{p{1.8cm}}{\raggedleft Root/Final upper bound} & 
 \multicolumn{1}{p{0.8cm}}{\raggedleft Best Result} & 
 \multicolumn{1}{p{0.7cm}}{\raggedleft \multirow{2}{*}{\centering Gap}} & 
 \multicolumn{1}{p{0.6cm}}{\raggedleft Run time(s)} \\
\noalign{\smallskip}\hline\noalign{\smallskip}
\textbf{C-1} & 194.72/194 & 194  & 0.00\% & 1 & 194.71/194 & 194 & 0.00\% & 0 \\
\textbf{C-2} & 600.62/599 & 599 & 0.00\% & 8 & 600.61/599 & 599 & 0.00\% & 1 \\
\textbf{C-3} & 643.35/629 & 629 & 0.00\% & 68 & 631.88/629 & 629 & 0.00\% & 3 \\
\textbf{C-4} & 334.01/318 & 318 & 0.00\% & 375 & 334.01/318 & 318 & 0.00\% & 39 \\
\textbf{C-5} & 903.61/894 & 881 & 1.45\% & -- & 899.58/888 & 887 & 0.11\% & -- \\
\textbf{C-6} & 951.06/939 & 928 & 1.17\% & -- & 946.05/928 & 928 & 0.00\% & 509 \\
\textbf{C-7} & 1211.39/1199 & 1132 & 5.59\% & -- & 1218.36/1198 & 1174 & 2.00\% & -- \\
\textbf{C-8} & 1364.65/1351 & 1300 & 3.77\% & -- & 1374.62/1348 & 1329 & 1.41\% & -- \\
\textbf{C-9} & 504.71/499 & 492 & 1.40\% & -- & 504.71/495 & 495 & 0.00\% & 912 \\
\textbf{C-10} & 1048.39/1042 & 1024 & 1.73\% & -- & 1055.39/1042 & 1038 & 0.38\% & -- \\
\textbf{C-11} & 1243.76/1239 & 1162 & 6.21\% & -- & 1243.76/1226 & 1219 & 0.57\% & -- \\
\textbf{C-12} & 1799.62/1793 & 1148 & 35.97\% & -- & 1806.59/1787 & 1735 & 2.91\% & -- \\
\textbf{C-13} & 2050.11/2050 & 1353 & 34.00\% & -- & 2060.03/2034 & 1974 & \textbf{2.95\%} & -- \\
\noalign{\smallskip}
\textbf{M-1} & 275.35/269 & 269 & 0.00\% & 5 & 272.04/269 & 269 & 0.00\% & 0 \\
\textbf{M-2} & 586.78/578 & 578 & 0.00\% & 60 & 590.40/578 & 578 & 0.00\% & 9 \\
\textbf{M-3} & 745.00/728 & 728 & 0.00\% & 496 & 739.31/728 & 728 & 0.00\% & 71 \\
\textbf{M-4} & 386.65/375 & 375 & 0.00\% & 51 & 382.50/375 & 375 & 0.00\% & 18 \\
\textbf{M-5} & 831.44/818 & 783 & 4.28\% &  & 823.26/807 & 800 & 0.87\% & -- \\
\textbf{M-6} & 1031.35/1018 & 982 & 3.54\% & -- & 1026.34/1009 & 1001 & 0.79\% & -- \\
\textbf{M-7} & 1558.55/1536 & 1510 & 1.69\% & -- & 1554.54/1533 & 1523 & 0.65\% & -- \\
\textbf{M-8} & 2037.25/2017 & 1964 & 2.63\% & -- & 2033.24/2008 & 1995 & 0.65\% & -- \\
\textbf{M-9} & 550.25/541 & 540 & 0.18\% & -- & 550.24/540 & 540 & 0.00\% & 49 \\
\textbf{M-10} & 1042.59/1041 & 969 & 6.92\% & -- & 1042.59/1036 & 1024 & 1.16\% & -- \\
\textbf{M-11} & 1371.50/1367 & 1059 & 22.53\% & -- & 1371.49/1360 & 1337 & 1.69\% & -- \\
\textbf{M-12} & 1981.15/1977 & 1627 & 17.70\% & -- & 1989.14/1970 & 1944 & 1.32\% & -- \\
\textbf{M-13} & 2552.34/2552 & 1746 & 31.58\% & -- & 2560.33/2544 & 2475 & 2.71\% & -- \\
\noalign{\smallskip}
\textbf{R-1} & 1344.55/1329 & 1329 & 0.00\% & 263 & 1332.42/1329 & 1329 & 0.00\% & 10 \\
\textbf{R-2} & 1816.51/1806 & 1806 & 0.00\% & 230 & 1809.70/1806 & 1806 & 0.00\% & 23 \\
\textbf{R-3} & 2366.41/2355 & 2355 & 0.00\% & 671 & 2358.70/2355 & 2355 & 0.00\% & 8 \\
\textbf{R-4} & 1550.02/1538 & 1520 & 1.17\% & -- & 1544.75/1528 & 1520 & 0.52\% & -- \\
\textbf{R-5} & 2039.12/2025 & 2006 & 0.94\% & -- & 2028.58/2015 & 2006 & 0.45\% & -- \\
\textbf{R-6} & 2581.86/2572 & 2552 & 0.78\% & -- & 2581.58/2568 & 2559 & 0.35\% & -- \\
\textbf{R-7} & 3171.19/3168 & 2835 & 10.51\% & -- & 3170.75/3157 & 3142 & 0.48\% & -- \\
\textbf{R-8} & 3635.33/3635 & 2984 & 17.91\% & -- & 3639.34/3621 & 3595 & 0.72\% & -- \\
\textbf{R-9} & 4217.33/4215 & 3533 & 16.18\% & -- & 4223.28/4204 & 4177 & 0.64\% & -- \\
\textbf{R-10} & 4700.81/4700 & 3984 & 15.23\% & -- & 4708.87/4688 & 4659 & 0.62\% & -- \\
\textbf{R-11} & -/-- & --  & -- & -- & 5189.80/5165 & 5125 & 0.77\% & -- \\
\noalign{\smallskip}\hline
\end{tabular}
\end{table}

These results can be improved further if we let Gurobi run for more than 6~hours. For example, for instance M-11, we can obtain the optimal solution with~205 assigned missions in 41,924~s, and a better solution with a total weight for assigned missions of~1340 in 35,523~s with a gap of~1.47\%. 

To directly compare instances, consider instances~C-13 and~M-13, which have the same available resources and the same number of missions for a scheduling horizon of 2~days. Both in the MILP and the improved MILP, the potential assignment opportunities $\paon$ and $\paot$ of missions for instance C-13 are higher than for instance M-13. This reflects that the improved model is more flexible in assigning a resource and an observation time to a mission. The size of the model scale of instance C-13 is slightly larger than that of instance M-13. Furthermore, for the improved MILP, the conflict indicator $\conf$ of instance C-13 is larger than that of instance M-13, meaning that the size of the model scale of the instance C-13 is much bigger than that of the instance M-13. 

\begin{figure}[!ht]
\centering
\includegraphics[width=119mm]{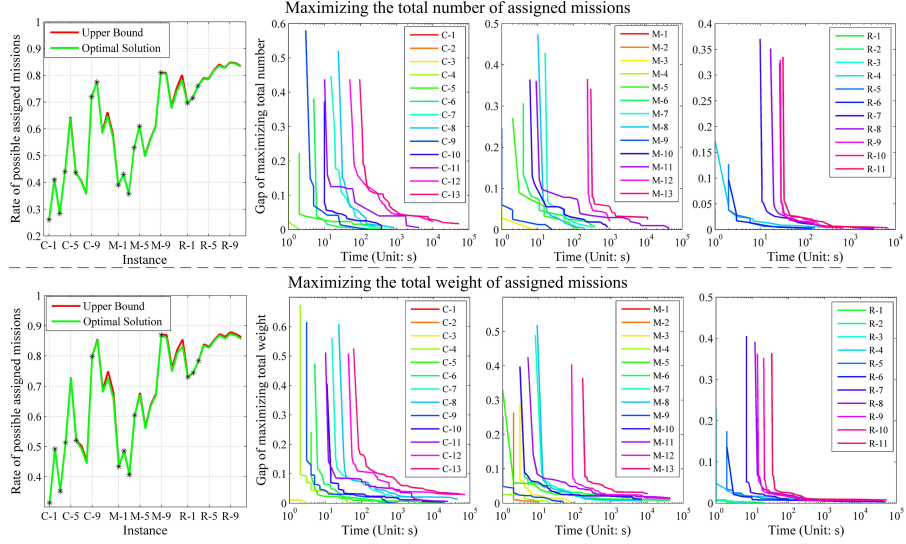}
\caption{Optimization results and the variation of computational efficiency}
\label{fig:ComEfficiency}
\end{figure}

In Fig.~\ref{fig:ComEfficiency}, we show further details about the problems and the optimization process. For instances with different computational complexity and different contention conflict for resources, the first column of the graph shows the potential capacity of resources for providing an upper bound as well as the performance of the proposed model for providing an optimal solution. The red line illustrates the maximum accomplishment rate of missions and the green line shows the best results that Gurobi achieves. The following three columns on the right of Fig.~\ref{fig:ComEfficiency} show the gap of the best solution obtained so far versus the runtime, and indicate that good upper bounds can be obtained fairly quickly. 

\section{Conclusions}\label{Sec:Conclusion}
In this paper, we addressed the problem of multi-satellite scheduling with limited observation capacities, which is one of the core problems to be solved for the effective utilization of the resources of satellite constellations. The key component of our approach is the detailed analysis of possible conflicts, which leads to stronger constraints in the MILP that significantly speed up the solution process. We find that the modelling of the problem with 5-index variables (thereby avoiding the standard ``Big-M'' approach) is more suitable for solving real-world instances in which most variables and constraints are not necessary for the model. 

In summary, for maximizing both the total number or the total weight of assigned missions, our experiments assess the correctness and effectiveness of the proposed MILP on several classes of problem instances. Good feasible solutions are obtained very fast, while the computation of true optimum solutions is also possible. The model thus provides a solid basis for designing a decision support system for scheduling satellite resources with imaging and communication tasks. 

\section*{Acknowledgments}
\noindent
We would like to thank the anonymous reviewers for their constructive, diligent, and detailed feedback that substantially improved the paper.

\noindent
This work is supported by National Natural Science Foundation of China under Grant No.41571403, the 13th Five-year Pre-research Project of Civil Aerospace in China, Joint Funds of Equipment Pre-Research and Ministry of Education of China under Grant No.6141A02022320, and Fundamental Research Funds for the Central Universities under Grant No.CUG2017G01 and No.CUG160207.

\section*{References}

\bibliography{mybibfile}


\begin{appendices}
\section{Generation of the Feasible Time Interval}\label{App:FTI}
\begin{figure}[!ht]
  \centering
  \includegraphics[width=86mm]{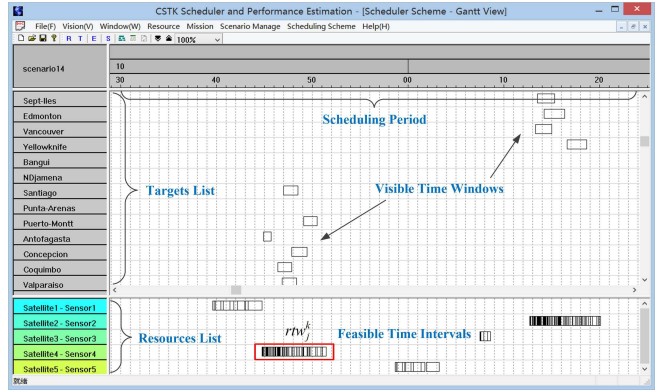}
  \caption{ Generation and distribution of feasible time intervals.}
  \label{fig:FeasibleTimeInterval1}
\end{figure}

If there is a free available time window, i.e., a sub-interval of a visible time window that does not overlap with any other windows, we assign it to the corresponding mission in the preprocessing step. The remaining visible time windows of the assigned mission are ignored. Therefore, all visible time windows that remain in the model are thus overlapping with at least one other visible time window. By combining all the overlapping visible time windows of missions on the same resource, the \emph{feasible time interval} is generated accordingly. The distributions of visible time windows and feasible time intervals are visualized in Fig.~\ref{fig:FeasibleTimeInterval1}. 

\begin{figure}[!ht]
  \centering
  \includegraphics[width=89mm]{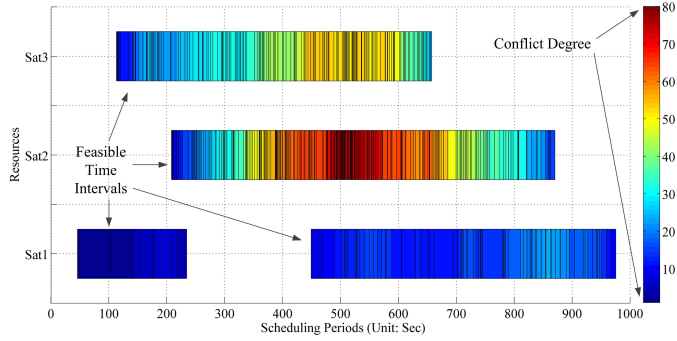}
  \caption{ Contention conflict distribution of feasible time intervals.}
  \label{fig:FeasibleTimeInterval2}
\end{figure}

Every resource~$R_j$ has several disjoint feasible time intervals, each consisting of pieces with a different number of overlapping visible time windows (\emph{the conflict degree}). Fig.~\ref{fig:FeasibleTimeInterval2} depicts the distribution of the conflict degree for the segments over the \emph{feasible time intervals}. Different colors denote the different conflict degrees of segments. 

\section{Calculation of Effective Feasible Time Subinterval}\label{App:SRTW}
Since we have knowledge of all corresponding missions that can be assigned to~$\srtw_{j}^{kl}$, we can calculate the value of each~$\srn_{j}^{kl}$. If all missions have the same observation duration time~$\dur$, then it can be computed as~$\srn_{j}^{kl} = \Bigl\lfloor\frac{|\srtw_{j}^{kl}|+\delta_j}{\dur+\delta_j}\Bigr\rfloor$. Otherwise, we iteratively assign a mission with the shortest observation duration time by taking the setup time constraint into account until it exceeds the capacity of~$\srtw_{j}^{kl}$ to obtain the value of~$\srn_{j}^{kl}$.

\begin{figure}[ht]
  \centering
  \includegraphics[width=119mm]{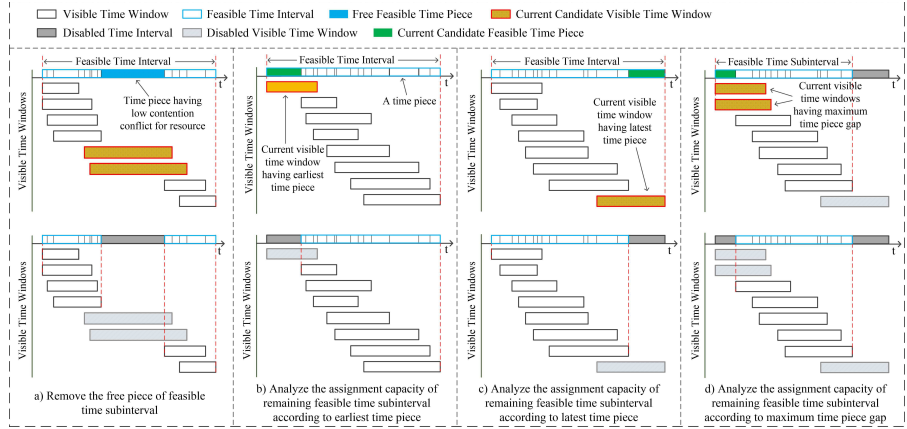}
  \caption{ Generation of effective feasible time subintervals.}
  \label{fig:FeasibleTimeInterval3}
\end{figure}

The generation of the \emph{effective feasible time subinterval} is handled as follows:
\begin{itemize}
\item 
Since the value of~$\srn_{j}^{kl}$ has already been calculated, we can perform partial assignments in a preprocessing phase. We remove the time-piece with the lowest conflict degree such that the number of the corresponding missions on~$\srtw_{j}^{kl}$ is less than or equal to~$\srn_{j}^{kl}$. All corresponding missions for this interval can be assigned directly, and the search space is decreased~(see Fig.~\ref{fig:FeasibleTimeInterval3} (a)). We denote the remaining time subintervals as \emph{effective feasible time subinterval}.

\item 
Considering the distribution of conflict degrees for feasible time intervals, we include more \emph{effective feasible time subintervals} and corresponding inequalities. To this end, we propose three operations. We iteratively ``remove'' a time-piece and its corresponding visible time windows according to the earliest start time, the latest end time, and the time that corresponds to the largest interval between the earliest start time and the latest end time separately~(see Fig.~\ref{fig:FeasibleTimeInterval3} b, c, and d).
\end{itemize}

\section{Instance Analysis}\label{App:InstAnaly}
In Table~\ref{tab:ResourceUtilisation}, we show the complexity of instances and the importance of each resource based on the utilization in different types of instances with differing conflict degree and distribution. In Table~\ref{tab:ResourceUtilisation}, $\delta$ denotes the maximum setup time of the resource, and $N$ is the total number of visible time windows available for each resource. The \emph{total visible time}~$T$ is the total visibility time over all visible time windows for a resource. The \emph{feasible observation time}~$F$ is the total time over feasible time intervals for a resource that can be assigned to missions (Due to the overlaps between visible time windows for the same resource, it is different from the \emph{total visible time}. Instead, it corresponds to the union of the overlapping visible time windows of a resource). 

\begin{table}[!htbp]
\scriptsize 
\caption{Resource Utilization}
\renewcommand{\tabcolsep}{8pt}
\label{tab:ResourceUtilisation}
\centering
\linespread{1.0}\selectfont
\begin{tabular}{lrrrrrrrrrr}
\hline\noalign{\smallskip}
\multirow{1}{*}{\centering Ins.} & 
\multirow{1}{*}{\centering Sat.} & 
\multirow{1}{*}{\centering Res.} & 
$\delta(s)$ & 
$N$ & 
$T$(s) & 
$F$(s) & 
$\rn$ & 
$\conf$ & 
$\paon$ & 
$\paot$ \\
\noalign{\smallskip}\hline\noalign{\smallskip}
\multirow{3}{*}{C-1} & HJ-1A & HIS & 30 & 44 & 5401.56 & 309.74 & 11 & 16.44 & \multirow{3}{*}{1.00} & \multirow{3}{*}{99.07} \\
 & HJ-1B & IRS & 30 & 45 & 4347.99 & 370.57 & 12 & 10.73 &  & \\
 & HJ-1C & SAR2 & 30 & 11 & 157.16 & 76.81 & 3 & 1.05 &  & \\
\noalign{\smallskip}
\multirow{3}{*}{C-2} & HJ-1A & HSI & 30 & 114 & 11692.25 & 1048.72 & 33 & 10.15 & \multirow{3}{*}{1.13} & \multirow{3}{*}{96.86} \\
 & HJ-1B & IRS & 30 & 83 & 7225.03 & 1219.97 & 38 & 4.92 &  & \\
 & HJ-1C & SAR2 & 30 & 28 & 455.16 & 312.85 & 17 & 0.45 &  & \\
\noalign{\smallskip}
\multirow{3}{*}{C-3} & HJ-1A & HSI & 30 & 177 & 16966.54 & 1460.07 & 42 & 10.62 & \multirow{3}{*}{1.21} & \multirow{3}{*}{99.86} \\
 & HJ-1B & IRS & 30 & 153 & 12466.22 & 1563.84 & 46 & 6.97 &  & \\
 & HJ-1C & SAR2 & 30 & 34 & 526.18 & 209.35 & 9 & 1.51 &  & \\
\noalign{\smallskip}
\multirow{5}{*}{C-4} & \multirow{2}{*}{HJ-1A} & CCD1 & 40 & 51 & 5056.82 & 400.81 & 9 & 11.62 & \multirow{5}{*}{1.95} & \multirow{5}{*}{191.69} \\
 &  & HSI & 30 & 44 & 5401.56 & 309.74 & 11 & 16.44 &  & \\
 & \multirow{2}{*}{HJ-1B} & CCD2 & 40 & 44 & 4205.78 & 359.60 & 9 & 10.70 &  & \\
 &  & IRS & 30 & 45 & 4347.99 & 370.57 & 12 & 10.73 &  & \\
 & HJ-1C & SAR2 & 30 & 11 & 157.16 & 76.81 & 3 & 1.05 &  & \\
\noalign{\smallskip}
\multirow{5}{*}{C-5} & \multirow{2}{*}{HJ-1A} & CCD1 & 40 & 136 & 12544.85 & 1408.84 & 34 & 7.90 & \multirow{5}{*}{2.31} & \multirow{5}{*}{202.06} \\
 &  & HSI & 30 & 114 & 11692.25 & 1048.72 & 33 & 10.15 &  & \\
 & \multirow{2}{*}{HJ-1B} & CCD2 & 40 & 101 & 8493.90 & 1302.05 & 33 & 5.52 &  & \\
 &  & IRS & 30 & 83 & 7225.03 & 1219.97 & 38 & 4.92 &  & \\
 & HJ-1C & SAR2 & 30 & 28 & 455.16 & 312.85 & 17 & 0.45 &  & \\
\noalign{\smallskip}
\multirow{5}{*}{C-6} & \multirow{2}{*}{HJ-1A} & CCD1 & 40 & 173 & 14221.61 & 1594.25 & 39 & 7.92 & \multirow{5}{*}{2.35} & \multirow{5}{*}{192.51} \\
 &  & HSI & 30 & 177 & 16966.54 & 1460.07 & 42 & 10.62 &  & \\
 & \multirow{2}{*}{HJ-1B} & CCD2 & 40 & 169 & 13573.70 & 1481.07 & 37 & 8.16 &  & \\
 &  & IRS & 30 & 153 & 12466.22 & 1563.84 & 46 & 6.97 &  & \\
 & HJ-1C & SAR2 & 30 & 34 & 526.18 & 209.35 & 9 & 1.51 &  & \\
\noalign{\smallskip}
\multirow{5}{*}{C-7} & \multirow{2}{*}{HJ-1A} & CCD1 & 40 & 251 & 20517.75 & 1783.73 & 36 & 10.50 & \multirow{5}{*}{2.56} & \multirow{5}{*}{221.70} \\
 &  & HSI & 30 & 251 & 25557.12 & 1930.87 & 45 & 12.24 &  & \\
 & \multirow{2}{*}{HJ-1B} & CCD2 & 40 & 237 & 20457.45 & 1785.02 & 43 & 10.46 &  & \\
 &  & IRS & 30 & 242 & 21479.04 & 1868.40 & 57 & 10.50 &  & \\
 & HJ-1C & SAR2 & 30 & 44 & 669.25 & 337.05 & 17 & 0.99 &  & \\
\noalign{\smallskip}
\multirow{5}{*}{C-8} & \multirow{2}{*}{HJ-1A} & CCD1 & 40 & 287 & 24015.23 & 1686.79 & 37 & 13.24 & \multirow{5}{*}{2.59} & \multirow{5}{*}{225.01} \\
 &  & HSI & 30 & 302 & 30849.21 & 1881.88 & 50 & 15.39 &  & \\
 & \multirow{2}{*}{HJ-1B} & CCD2 & 40 & 323 & 27687.92 & 2027.46 & 49 & 12.66 &  & \\
 &  & IRS & 30 & 330 & 29139.61 & 2110.58 & 62 & 12.81 &  & \\
 & HJ-1C & SAR2 & 30 & 53 & 811.56 & 364.61 & 18 & 1.23 &  & \\
\noalign{\smallskip}
\multirow{5}{*}{C-9} & \multirow{2}{*}{HJ-1A} & CCD1 & 40 & 125 & 11638.60 & 813.02 & 19 & 13.32 & \multirow{5}{*}{3.67} & \multirow{5}{*}{366.30} \\
 &  & HSI & 30 & 88 & 10659.17 & 623.03 & 21 & 16.11 &  & \\
 & \multirow{2}{*}{HJ-1B} & CCD2 & 40 & 71 & 6972.31 & 548.55 & 14 & 11.71 &  & \\
 &  & IRS & 30 & 72 & 7202.65 & 560.91 & 18 & 11.84 &  & \\
 & HJ-1C & SAR2 & 30 & 11 & 157.16 & 76.81 & 3 & 1.05 &  & \\
\noalign{\smallskip}
\multirow{5}{*}{C-10} & \multirow{2}{*}{HJ-1A} & CCD1 & 40 & 253 & 22024.57 & 2223.94 & 55 & 8.90 & \multirow{5}{*}{3.40} & \multirow{5}{*}{301.76} \\
 &  & HSI & 30 & 181 & 18953.14 & 1593.73 & 49 & 10.89 &  & \\
 & \multirow{2}{*}{HJ-1B} & CCD2 & 40 & 123 & 10471.08 & 1724.89 & 43 & 5.07 &  & \\
 &  & IRS & 30 & 95 & 8448.97 & 1463.36 & 44 & 4.77 &  & \\
 & HJ-1C & SAR2 & 30 & 28 & 455.16 & 312.85 & 17 & 0.45 &  & \\
\noalign{\smallskip}
\multirow{5}{*}{C-11} & \multirow{2}{*}{HJ-1A} & CCD1 & 40 & 369 & 29721.34 & 2673.48 & 65 & 10.12 & \multirow{5}{*}{3.69} & \multirow{5}{*}{310.26} \\
 &  & HSI & 30 & 273 & 26917.49 & 1899.59 & 56 & 13.17 &  & \\
 & \multirow{2}{*}{HJ-1B} & CCD2 & 40 & 248 & 20545.92 & 2174.95 & 55 & 8.45 &  & \\
 &  & IRS & 30 & 182 & 15365.81 & 1851.21 & 55 & 7.30 &  & \\
 & HJ-1C & SAR2 & 30 & 34 & 526.18 & 209.35 & 9 & 1.51 &  & \\
\noalign{\smallskip}
\multirow{5}{*}{C-12} & \multirow{2}{*}{HJ-1A} & CCD1 & 40 & 543 & 44443.30 & 3975.35 & 86 & 10.18 & \multirow{5}{*}{4.34} & \multirow{5}{*}{377.97} \\
 &  & HSI & 30 & 378 & 38680.16 & 3403.81 & 81 & 10.36 &  & \\
 & \multirow{2}{*}{HJ-1B} & CCD2 & 40 & 451 & 38735.80 & 3582.40 & 79 & 9.81 &  & \\
 &  & IRS & 30 & 321 & 28658.97 & 2717.88 & 80 & 9.54 &  & \\
 & HJ-1C & SAR2 & 30 & 44 & 669.25 & 337.05 & 17 & 0.99 &  & \\
\noalign{\smallskip}
\multirow{5}{*}{C-13} & \multirow{2}{*}{HJ-1A} & CCD1 & 40 & 677 & 56901.40 & 3852.66 & 88 & 13.77 & \multirow{5}{*}{4.31} & \multirow{5}{*}{376.24} \\
 &  & HSI & 30 & 447 & 45886.18 & 3407.96 & 91 & 12.46 &  & \\
 & \multirow{2}{*}{HJ-1B} & CCD2 & 40 & 569 & 48600.95 & 3898.63 & 89 & 11.47 &  & \\
 &  & IRS & 30 & 407 & 35921.75 & 3062.02 & 89 & 10.73 &  & \\
 & HJ-1C & SAR2 & 30 & 53 & 811.56 & 364.61 & 18 & 1.23 &  & \\
\noalign{\smallskip}\hline
\end{tabular}
\end{table}

\begin{table}[!htbp]
\scriptsize 
\centering
\renewcommand{\tabcolsep}{8pt}
\linespread{1.0}\selectfont
\begin{tabular}{lrrrrrrrrrr}
\hline\noalign{\smallskip}
\multirow{1}{*}{\centering Ins.} & 
\multirow{1}{*}{\centering Sat.} & 
\multirow{1}{*}{\centering Res.} & 
$\delta(s)$ & 
$N$ & 
$T$(s) & 
$F$(s) & 
$\rn$ & 
$\conf$ & 
$\paon$ & 
$\paot$ \\
\noalign{\smallskip}\hline\noalign{\smallskip}
\multirow{3}{*}{M-1} & HJ-1A & HIS & 30 & 45 & 5138.73 & 567.13 & 14 & 8.06 & \multirow{3}{*}{1.07} & \multirow{3}{*}
{101.40} \\
 & HJ-1B & IRS & 30 & 54 & 4874.60 & 1595.55 & 29 & 2.06 &  & \\
 & HJ-1C & SAR2 & 30 & 8 & 126.23 & 67.39 & 3 & 0.87 &  & \\
\noalign{\smallskip}
\multirow{3}{*}{M-2} & HJ-1A & HSI & 30 & 93 & 9381.23 & 1460.33 & 37 & 5.42 & \multirow{3}{*}{1.16} & \multirow{3}{*}
{100.17} \\
 & HJ-1B & IRS & 30 & 120 & 10356.08 & 2194.82 & 59 & 3.72 &  & \\
 & HJ-1C & SAR2 & 30 & 18 & 296.09 & 140.49 & 6 & 1.11 &  & \\
\noalign{\smallskip}
\multirow{3}{*}{M-3} & HJ-1A & HSI & 30 & 151 & 15150.84 & 1738.97 & 44 & 7.71 & \multirow{3}{*}{1.19} & \multirow{3}{*}
{103.84} \\
 & HJ-1B & IRS & 30 & 184 & 15645.23 & 3223.84 & 77 & 3.85 &  & \\
 & HJ-1C & SAR2 & 30 & 22 & 355.19 & 148.52 & 6 & 1.39 &  & \\
\noalign{\smallskip}
\multirow{5}{*}{M-4} & \multirow{2}{*}{HJ-1A} & CCD1 & 40 & 43 & 4283.88 & 406.60 & 9 & 9.54 & \multirow{5}{*}{2.07} & 
\multirow{5}{*}{194.17} \\
 &  & HSI & 30 & 45 & 5138.73 & 567.13 & 14 & 8.06 &  & \\
 & \multirow{2}{*}{HJ-1B} & CCD2 & 40 & 57 & 4994.01 & 1503.96 & 27 & 2.32 &  & \\
 &  & IRS & 30 & 54 & 4874.60 & 1595.55 & 29 & 2.06 &  & \\
 & HJ-1C & SAR2 & 30 & 8 & 126.23 & 67.39 & 3 & 0.87 &  & \\
\noalign{\smallskip}
\multirow{5}{*}{M-5} & \multirow{2}{*}{HJ-1A} & CCD1 & 40 & 97 & 8688.12 & 1440.46 & 34 & 5.03 & \multirow{5}{*}{2.31} & 
\multirow{5}{*}{199.85} \\
 &  & HSI & 30 & 93 & 9381.23 & 1460.33 & 37 & 5.42 &  & \\
 & \multirow{2}{*}{HJ-1B} & CCD2 & 40 & 134 & 11248.59 & 2111.44 & 51 & 4.33 &  & \\
 &  & IRS & 30 & 120 & 10356.08 & 2194.82 & 59 & 3.72 &  & \\
 & HJ-1C & SAR2 & 30 & 18 & 296.09 & 140.49 & 6 & 1.11 &  & \\
\noalign{\smallskip}
\multirow{5}{*}{M-6} & \multirow{2}{*}{HJ-1A} & CCD1 & 40 & 160 & 13827.02 & 1656.37 & 39 & 7.35 & \multirow{5}{*}{2.39} & 
\multirow{5}{*}{205.43} \\
 &  & HSI & 30 & 151 & 15150.84 & 1738.97 & 44 & 7.71 &  & \\
 & \multirow{2}{*}{HJ-1B} & CCD2 & 40 & 199 & 16650.92 & 3098.04 & 66 & 4.37 &  & \\
 &  & IRS & 30 & 184 & 15645.23 & 3223.84 & 77 & 3.85 &  & \\
 & HJ-1C & SAR2 & 30 & 22 & 355.19 & 148.52 & 6 & 1.39 &  & \\
\noalign{\smallskip}
\multirow{5}{*}{M-7} & \multirow{2}{*}{HJ-1A} & CCD1 & 40 & 295 & 24607.41 & 5392.2 & 96 & 3.56 & \multirow{5}{*}{2.30} & \multirow{5}{*}{194.53} \\
 &  & HSI & 30 & 291 & 28775.26 & 2633.39 & 64 & 9.93 &  & \\
 & \multirow{2}{*}{HJ-1B} & CCD2 & 40 & 210 & 17288.81 & 4265.57 & 74 & 3.05 &  & \\
 &  & IRS & 30 & 81 & 6485.44 & 3207.25 & 57 & 1.02 &  & \\
 & HJ-1C & SAR2 & 30 & 41 & 655.26 & 485.40 & 25 & 0.35 &  & \\
\noalign{\smallskip}
\multirow{5}{*}{M-8} & \multirow{2}{*}{HJ-1A} & CCD1 & 40 & 334 & 27847.32 & 6509.20 & 122 & 3.28 & \multirow{5}{*}{2.16} & \multirow{5}{*}{178.41} \\
 &  & HSI & 30 & 300 & 28891.38 & 3416.69 & 83 & 7.46 &  & \\
 & \multirow{2}{*}{HJ-1B} & CCD2 & 40 & 308 & 24596.35 & 6820.68 & 126 & 2.61 &  & \\
 &  & IRS & 30 & 93 & 7146.66 & 3298.15 & 62 & 1.17 &  & \\
 & HJ-1C & SAR2 & 30 & 45 & 721.53 & 573.52 & 30 & 0.26 &  & \\
\noalign{\smallskip}
\multirow{5}{*}{M-9} & \multirow{2}{*}{HJ-1A} & CCD1 & 40 & 121 & 11012.80 & 776.74 & 19 & 13.18 & \multirow{5}{*}{3.93} & 
\multirow{5}{*}{374.05} \\
 &  & HSI & 30 & 90 & 10397.46 & 982.28 & 27 & 9.59 &  & \\
 & \multirow{2}{*}{HJ-1B} & CCD2 & 40 & 90 & 8094.21 & 2825.45 & 47 & 1.86 &  & \\
 &  & IRS & 30 & 84 & 7774.62 & 2883.95 & 48 & 1.70 &  & \\
 & HJ-1C & SAR2 & 30 & 8 & 126.23 & 67.39 & 3 & 0.87 &  & \\
\noalign{\smallskip}
\multirow{5}{*}{M-10} & \multirow{2}{*}{HJ-1A} & CCD1 & 40 & 281 & 23817.80 & 3359.23 & 80 & 6.09 & \multirow{5}{*}{4.07} & 
\multirow{5}{*}{356.65} \\
 &  & HSI & 30 & 170 & 17551.99 & 2243.98 & 56 & 6.82 &  & \\
 & \multirow{2}{*}{HJ-1B} & CCD2 & 40 & 201 & 17003.69 & 3642.71 & 85 & 3.67 &  & \\
 &  & IRS & 30 & 144 & 12660.71 & 2902.44 & 74 & 3.36 &  & \\
 & HJ-1C & SAR2 & 30 & 18 & 296.09 & 140.49 & 6 & 1.11 &  & \\
\noalign{\smallskip}
\multirow{5}{*}{M-11} & \multirow{2}{*}{HJ-1A} & CCD1 & 40 & 407 & 33683.82 & 3623.79 & 87 & 8.30 & \multirow{5}{*}{4.17} & 
\multirow{5}{*}{361.57} \\
 &  & HSI & 30 & 255 & 25859.69 & 2577.14 & 65 & 9.03 &  & \\
 & \multirow{2}{*}{HJ-1B} & CCD2 & 40 & 332 & 27988.76 & 5624.21 & 114 & 3.98 &  & \\
 &  & IRS & 30 & 236 & 20582.38 & 4786.82 & 103 & 3.30 &  & \\
 & HJ-1C & SAR2 & 30 & 22 & 355.19 & 148.52 & 6 & 1.39 &  & \\
\noalign{\smallskip}
\multirow{5}{*}{M-12} & \multirow{2}{*}{HJ-1A} & CCD1 & 40 & 524 & 43726.61 & 9351.85 & 165 & 3.68 & \multirow{5}{*}{3.99} & \multirow{5}{*}{339.35} \\
 &  & HSI & 30 & 479 & 47108.39 & 5264.96 & 122 & 7.95 &  & \\
 & \multirow{2}{*}{HJ-1B} & CCD2 & 40 & 423 & 34216.82 & 7741.70 & 142 & 3.42 &  & \\
 &  & IRS & 30 & 127 & 10015.32 & 5094.58 & 94 & 0.97 &  & \\
 & HJ-1C & SAR2 & 30 & 42 & 673.92 & 504.06 & 26 & 0.34 &  & \\
\noalign{\smallskip}
\multirow{5}{*}{M-13} & \multirow{2}{*}{HJ-1A} & CCD1 & 40 & 667 & 56629.50 & 14082.51 & 258 & 3.02 & \multirow{5}{*}{4.07} & \multirow{5}{*}{347.95} \\
 &  & HSI & 30 & 576 & 56857.17 & 9312.37 & 206 & 5.11 &  & \\
 & \multirow{2}{*}{HJ-1B} & CCD2 & 40 & 575 & 46201.02 & 13102.91 & 238 & 2.53 &  & \\
 &  & IRS & 30 & 172 & 13566.61 & 7074.17 & 128 & 0.92 &  & \\
 & HJ-1C & SAR2 & 30 & 45 & 721.53 & 573.52 & 30 & 0.26 &  & \\
\noalign{\smallskip}\hline
\end{tabular}
\end{table}

\begin{table}[!htbp]
\scriptsize 
\centering
\renewcommand{\tabcolsep}{8pt}
\linespread{1.0}\selectfont
\begin{tabular}{lrrrrrrrrrr}
\hline\noalign{\smallskip}
\multirow{1}{*}{\centering Ins.} & 
\multirow{1}{*}{\centering Sat.} & 
\multirow{1}{*}{\centering Res.} & 
$\delta(s)$ & 
$N$ & 
$T$(s) & 
$F$(s) & 
$\rn$ & 
$\conf$ & 
$\paon$ & 
$\paot$ \\
\noalign{\smallskip}\hline\noalign{\smallskip}
\multirow{4}{*}{R-1} & HJ-1A & HSI & 30 & 170 & 17219.78 & 7210.43 & 118 & 1.39 & \multirow{4}{*}{1.31} & \multirow{4}{*}{107.36} \\
 & HJ-1B & IRS & 30 & 166 & 14165.54 & 5382.46 & 94 & 1.63 &  & \\
 & \multirow{2}{*}{HJ-1C} & SAR1 & 25 & 11 & 116.99 & 116.99 & 11 & 0.00 &  & \\
 &  & SAR2 & 30 & 46 & 706.49 & 653.70 & 39 & 0.08 &  & \\
\noalign{\smallskip}
\multirow{4}{*}{R-2} & HJ-1A & HSI & 30 & 210 & 21363.18 & 8508.24 & 142 & 1.51 & \multirow{4}{*}{1.28} & \multirow{4}{*}{105.77} \\
 & HJ-1B & IRS & 30 & 231 & 19933.08 & 9559.21 & 152 & 1.09 &  & \\
 & \multirow{2}{*}{HJ-1C} & SAR1 & 25 & 12 & 126.75 & 126.75 & 12 & 0.00 &  & \\
 &  & SAR2 & 30 & 57 & 883.27 & 777.86 & 47 & 0.14 &  & \\
\noalign{\smallskip}
\multirow{4}{*}{R-3} & HJ-1A & HSI & 30 & 247 & 24748.59 & 10327.56 & 179 & 1.40 & \multirow{4}{*}{1.26} & \multirow{4}{*}{101.96} \\
 & HJ-1B & IRS & 30 & 285 & 24789.20 & 12415.51 & 205 & 1.00 &  & \\
 & \multirow{2}{*}{HJ-1C} & SAR1 & 25 & 20 & 211.59 & 204.68 & 18 & 0.03 &  & \\
 &  & SAR2 & 30 & 79 & 1233.10 & 1062.72 & 63 & 0.16 &  & \\
\noalign{\smallskip}
\multirow{6}{*}{R-4} & \multirow{2}{*}{HJ-1A} & CCD1 & 40 & 170 & 14257.38 & 6755.10 & 112 & 1.11 & \multirow{6}{*}{2.43} & \multirow{6}{*}{200.96} \\
 &  & HSI & 30 & 170 & 17219.78 & 7210.43 & 118 & 1.39 &  & \\
 & \multirow{2}{*}{HJ-1B} & CCD2 & 40 & 166 & 13822.88 & 5239.04 & 86 & 1.64 &  & \\
 &  & IRS & 30 & 166 & 14165.54 & 5382.46 & 94 & 1.63 &  & \\
 & \multirow{2}{*}{HJ-1C} & SAR1 & 25 & 11 & 116.99 & 116.99 & 11 & 0.00 &  & \\
 &  & SAR2 & 30 & 46 & 706.49 & 653.70 & 39 & 0.08 &  & \\
\noalign{\smallskip}
\multirow{6}{*}{R-5} & \multirow{2}{*}{HJ-1A} & CCD1 & 40 & 245 & 20460.47 & 10263.72 & 171 & 0.99 & \multirow{6}{*}{2.47} & \multirow{6}{*}{205.99} \\
 &  & HSI & 30 & 210 & 21363.18 & 8508.24 & 142 & 1.51 &  & \\
 & \multirow{2}{*}{HJ-1B} & CCD2 & 40 & 234 & 19627.92 & 9203.93 & 142 & 1.13 &  & \\
 &  & IRS & 30 & 231 & 19933.08 & 9559.21 & 152 & 1.09 &  & \\
 & \multirow{2}{*}{HJ-1C} & SAR1 & 25 & 12 & 126.75 & 126.75 & 12 & 0.00 &  & \\
 &  & SAR2 & 30 & 57 & 883.27 & 777.86 & 47 & 0.14 &  & \\
\noalign{\smallskip}
\multirow{6}{*}{R-6} & \multirow{2}{*}{HJ-1A} & CCD1 & 40 & 298 & 24794.29 & 12754.80 & 221 & 0.94 & \multirow{6}{*}{2.44} & \multirow{6}{*}{200.60} \\
 &  & HSI & 30 & 247 & 24748.59 & 10327.56 & 179 & 1.40 &  & \\
 & \multirow{2}{*}{HJ-1B} & CCD2 & 40 & 291 & 24523.34 & 12205.38 & 199 & 1.01 &  & \\
 &  & IRS & 30 & 285 & 24789.20 & 12415.51 & 205 & 1.00 &  & \\
 & \multirow{2}{*}{HJ-1C} & SAR1 & 25 & 20 & 211.59 & 204.68 & 18 & 0.03 &  & \\
 &  & SAR2 & 30 & 79 & 1233.10 & 1062.72 & 63 & 0.16 &  & \\
\noalign{\smallskip}
\multirow{6}{*}{R-7} & \multirow{2}{*}{HJ-1A} & CCD1 & 40 & 394 & 32805.30 & 15058.92 & 272 & 1.18 & \multirow{6}{*}{2.62} & \multirow{6}{*}{218.61} \\
 &  & HSI & 30 & 359 & 35980.50 & 13316.92 & 252 & 1.70 &  & \\
 & \multirow{2}{*}{HJ-1B} & CCD2 & 40 & 378 & 31809.02 & 14261.62 & 247 & 1.23 &  & \\
 &  & IRS & 30 & 332 & 28953.59 & 12883.02 & 224 & 1.25 &  & \\
 & \multirow{2}{*}{HJ-1C} & SAR1 & 25 & 21 & 220.32 & 213.41 & 19 & 0.03 &  & \\
 &  & SAR2 & 30 & 89 & 1394.67 & 1178.06 & 71 & 0.18 &  & \\
\noalign{\smallskip}
\multirow{6}{*}{R-8} & \multirow{2}{*}{HJ-1A} & CCD1 & 40 & 468 & 39119.41 & 17147.12 & 315 & 1.28 & \multirow{6}{*}{2.72} & \multirow{6}{*}{228.59} \\
 &  & HSI & 30 & 449 & 44694.12 & 15967.84 & 312 & 1.80 &  & \\ 
 & \multirow{2}{*}{HJ-1B} & CCD2 & 40 & 492 & 41659.75 & 17020.37 & 310 & 1.45 &  & \\
 &  & IRS & 30 & 374 & 32698.60 & 13032.66 & 229 & 1.51 &  & \\
 & \multirow{2}{*}{HJ-1C} & SAR1 & 25 & 23 & 239.07 & 232.16 & 21 & 0.03 &  & \\
 &  & SAR2 & 30 & 101 & 1602.63 & 1326.67 & 79 & 0.21 &  & \\
\noalign{\smallskip}
\multirow{6}{*}{R-9} & \multirow{2}{*}{HJ-1A} & CCD1 & 40 & 597 & 49603.80 & 20593.61 & 395 & 1.41 & \multirow{6}{*}{2.88} & \multirow{6}{*}{241.49} \\
 &  & HSI & 30 & 538 & 53486.70 & 19275.32 & 384 & 1.77 &  & \\
 & \multirow{2}{*}{HJ-1B} & CCD2 & 40 & 620 & 51972.71 & 20314.63 & 384 & 1.56 &  & \\
 &  & IRS & 30 & 416 & 36115.20 & 14705.63 & 265 & 1.46 &  & \\
 & \multirow{2}{*}{HJ-1C} & SAR1 & 25 & 25 & 258.94 & 252.03 & 23 & 0.03 &  & \\
 &  & SAR2 & 30 & 111 & 1758.41 & 1478.82 & 88 & 0.19 &  & \\
\noalign{\smallskip}
\multirow{6}{*}{R-10} & \multirow{2}{*}{HJ-1A} & CCD1 & 40 & 666 & 55439.42 & 23101.46 & 443 & 1.40 & \multirow{6}{*}{2.86} & \multirow{6}{*}{237.52} \\
 &  & HSI & 30 & 578 & 57447.53 & 20896.70 & 415 & 1.75 &  & \\
 & \multirow{2}{*}{HJ-1B} & CCD2 & 40 & 696 & 58559.76 & 22731.15 & 435 & 1.58 &  & \\
 &  & IRS & 30 & 459 & 39784.35 & 16460.32 & 303 & 1.42 &  & \\
 & \multirow{2}{*}{HJ-1C} & SAR1 & 25 & 30 & 307.84 & 300.93 & 28 & 0.02 &  & \\
 &  & SAR2 & 30 & 141 & 2230.07 & 1889.91 & 113 & 0.18 &  & \\
\noalign{\smallskip}
\multirow{6}{*}{R-11} & \multirow{2}{*}{HJ-1A} & CCD1 & 40 & 783 & 65837.40 & 24557.99 & 477 & 1.68 & \multirow{6}{*}{2.92} & \multirow{6}{*}{244.30} \\
 &  & HSI & 30 & 657 & 65689.49 & 22099.31 & 444 & 1.97 &  & \\
 & \multirow{2}{*}{HJ-1B} & CCD2 & 40 & 784 & 65688.80 & 23582.36 & 452 & 1.79 &  & \\
 &  & IRS & 30 & 511 & 44351.37 & 17214.34 & 320 & 1.58 &  & \\
 & \multirow{2}{*}{HJ-1C} & SAR1 & 25 & 33 & 343.53 & 336.61 & 31 & 0.02 &  & \\
 &  & SAR2 & 30 & 150 & 2385.59 & 2008.28 & 118 & 0.19 &  & \\
\noalign{\smallskip}\hline
\end{tabular}
\end{table}

\end{appendices}

\end{document}